\documentclass[final,leqno]{siamltex}
\usepackage{amsmath}
\usepackage{amssymb}
\usepackage{float}
\usepackage{graphicx}
\usepackage{textcomp}
\usepackage{subfig}
\usepackage{srcltx}
\usepackage{hyperref}
\hypersetup{draft=false, colorlinks=true, linkcolor=black, urlcolor=blue, citecolor=black}
\newtheorem{remark}[theorem]{Remark}
\newtheorem{example}[theorem]{Example}
\newtheorem{algorithm}[theorem]{Algorithm}
\newtheorem{assumption}[theorem]{Assumption}

\newcommand{\RR}{\mathbb{R}}

\newcommand{\NN}{\mathbb{N}}

\newcommand{\laplace}{\Delta}

\newcommand{\ctc}{c_{\mathrm{tc}}}

\newcommand{\xs}{x^*}

\newcommand{\Xs}{X^*}

\author{Anne Wald\thanks{Department of Mathematics, Saarland University, PO Box 15 11 50, 66041 Saarbr\"ucken, Germany ({\tt anne.wald@num.uni-sb.de}).}}

\title{A fast subspace optimization method for nonlinear inverse problems in Banach spaces with an application in parameter identification}
        
\begin{document}

\maketitle

\begin{abstract}
 We introduce and analyze a fast iterative method based on sequential Bregman projections for nonlinear inverse problems in Banach spaces. The key idea, in contrast to the standard Landweber method, is to use multiple search directions per iteration in combination with a regulation of the step width in order to reduce the total number of iterations. This method is suitable for both exact and noisy data. In the latter case, we obtain a regularization method. An algorithm with two search directions is used for the numerical identification of a parameter in an elliptic boundary value problem.
\end{abstract}



\section{Introduction}
Inverse problems, described by operator equations
\begin{displaymath}
 F(x) = y, \quad F: \mathcal{D}(F) \subseteq X \rightarrow Y,
\end{displaymath}
where $\mathcal{D}(F)$ denotes the domain of $F$, are usually classified according to the properties of the forward operator $F$. First of all, we distinguish between linear and nonlinear inverse problems. If $F$ is linear, we call the respective inverse problem linear. This class of problems has been discussed in a wide range of literature, see for example \cite{ehn00, louis89, skhk12}. The class of nonlinear inverse problems, where the forward operator $F$ is nonlinear, is addressed in, e.g., \cite{bh00, kns_itreg}. A further classification arises from the type of spaces in which an inverse problem is formulated: we differentiate between Hilbert and Banach spaces. In Hilbert spaces, we have a range of tools such as scalar products and orthogonal projections that play an important role in reconstruction and regularization methods. Also, we usually identify a Hilbert space with its dual space, which is a handy property for reconstruction techniques. The aforementioned references are mainly dealing with inverse problems in 
Hilbert spaces. If the inverse problem is defined in a Banach space setting, we have to find a solution without the help of these tools. An overview of Banach spaces techniques is given in \cite{skhk12}. \\
In this article, we want to extend the sequential subspace optimization (SESOP) methods, which were first developed to solve linear systems of equations in finite dimensional vector spaces \cite{gnmz05} and have been introduced for linear inverse problems in Hilbert and Banach spaces \cite{ss09, ssl08} as well as for nonlinear inverse problems in Hilbert spaces \cite{ws16, awts17}, to the Banach space setting. \\
The main purpose of SESOP methods is to reduce the number of necessary iterations in order to speed up the reconstruction by increasing the dimension of the search space, i.e., the (affine) subspace of $X$, in which the subsequent iterate is calculated. We thus use multiple search directions per iteration, but in contrast to the conjugate gradient (CG) method \cite{mhes52} or the generalized minimal residual (GMRES) method \cite{ysms86} we do not use Krylov spaces or increasing search space dimensions. Instead, we use search spaces with a bounded dimension. \\
As discussed in \cite{ss09, ssl08, ws16}, the SESOP methods allow a geometrical interpretation: The new iterate is calculated as the projection of the current iterate onto the intersection of hyperplanes or stripes that are determined by the search directions and whose width depends on the noise level in the data and on the nonlinearity of the forward operator. Hence, the character of the respective inverse problem is strongly reflected in the method: On the one hand, the (non)linearity of the forward operator defines the geometry of the stripes, whereas, on the other hand, the type of spaces determines the type of projection that is used: In a Hilbert space setting, we use orthogonal and metric projections. In a Banach space setting, however, we have to use Bregman projections. \\
In Section \ref{sec_prelim} we introduce the required notation and give an overview of the essential tools such as duality mappings, Bregman projections, and Bregman distances, along with some properties of nonlinear operators. Afterwards, in Section \ref{sec_sesop}, we introduce the SESOP and regularizing SESOP methods for nonlinear inverse problems in Banach spaces and show convergence and regularization properties for a certain selection of search directions. The method is then evaluated by solving a well-understood nonlinear parameter identification problem in Section \ref{sec_numex}. Finally, we discuss the findings of this article and give an outlook to possible extensions of this work.

\section{Mathematical setup} \label{sec_prelim}
Throughout this paper, let $(X,\lVert \cdot \rVert_X)$ and $(Y,\lVert \cdot \rVert_Y)$ denote real Banach spaces with their respective norms. Their respective duals are denoted by $(X^*, \lVert\cdot\rVert_{X^*})$ and $(Y^*, \lVert\cdot\rVert_{Y^*})$. The \emph{dual pairing} of $x\in X$ and $\xs \in \Xs$ is defined as
\begin{equation}
 \langle \xs, x \rangle_{\Xs \times X} := \langle x, \xs \rangle_{X \times \Xs} := \xs(x).
\end{equation}
For a better readability we drop the subscripts whenever confusion is not possible. 

\vspace*{1ex}

For $p>1$, the \emph{conjugate exponent} $p^*>1$ is determined by
\begin{displaymath}
 \frac{1}{p} + \frac{1}{p^*} = 1.
\end{displaymath}

\vspace*{1ex}

We consider the operator equation
\begin{equation} \label{op_eq}
 F(x) = y, \quad F: \mathcal{D}(F) \subseteq X \rightarrow Y,
\end{equation}
where the nonlinear operator $F$ is continuous and Fr\'echet differentiable in a ball $B_{\rho}(x_0) \subseteq \mathcal{D}(F)$ centered about $x_0 \in \mathcal{D}(F)$ with radius $\rho>0$. We postulate that $F$ fulfills the \emph{tangential cone condition}
\begin{equation}\label{tcc}
 \left\lVert F(x) - F(\tilde{x}) - F'(x)(x-\tilde{x}) \right\rVert \leq c_{\mathrm{tc}} \left\lVert F(x) - F(\tilde{x}) \right\rVert
\end{equation}
with $0 \leq \ctc < 1$ for $x,\tilde{x} \in B_{\rho}(x_0)$. Furthermore, let
\begin{displaymath}
 M_{F(x)=y} := \left\lbrace x \in X \, : \, F(x) = y \right\rbrace
\end{displaymath}
be the solution set of the operator equation \eqref{op_eq}. In case only noisy data $y^{\delta}$ is given, we assume that the noise level $\delta > 0$ fulfills
\begin{displaymath}
 \big\lVert y - y^{\delta} \big\rVert \leq \delta.
\end{displaymath}

\vspace*{2ex}

\begin{proposition}\label{lemma_prop_F}
 The validity of the tangential cone condition \eqref{tcc} in a ball $B_{\rho}(x_0)$ for some $x_0 \in X$ implies that the Fr\'echet derivative fulfills $F'(x) = 0$ for some $x \in B_{\rho}(x_0)$ if and only if $F$ is constant in $B_{\rho}(x_0)$.
\end{proposition}

\vspace*{2ex}

\begin{proof}
 Let $F'(x) = 0$ for some $x \in B_{\rho}(x_0)$. For all $\tilde{x} \in B_{\rho}(x_0)$ the tangential cone condition yields
 \begin{displaymath}
  \left\lVert F(x) - F(\tilde{x}) \right\rVert \leq c_{\mathrm{tc}} \left\lVert F(x) - F(\tilde{x}) \right\rVert,
 \end{displaymath}
 which can only be satisfied if $F(\tilde{x}) = F(x)$ in $B_{\rho}(x_0)$, i.e., $F$ is constant in $B_{\rho}(x_0)$. \\
 Now let $F$ be constant in $B_{\rho}(x_0)$. The tangential cone condition now yields
 \begin{displaymath}
  \left\lVert F'(x)(x-\tilde{x}) \right\rVert \leq 0
 \end{displaymath}
 for all $x,\tilde{x} \in B_{\rho}(x_0)$, implying $F'(x) = 0$ for $x \in B_{\rho}(x_0)$.
\end{proof}

\vspace*{2ex}

We further postulate that 
\begin{equation} \label{eq_c_F}
 \lVert F'(x) \rVert \leq c_F \quad \text{for all} \ x \in B_{\rho}(x_0)
\end{equation}
for a constant $c_F \geq 0$. According to Proposition \ref{lemma_prop_F}, the case $c_F=0$ is not interesting, so we use $c_F>0$.

Additionally, let $F$ be weakly sequentially closed in $B_{\rho}(x_0)$, i.e., for $x_n, x \in B_{\rho}(x_0)$, $n \in \NN$, with $x_n \rightharpoonup x$ and $F(x_n) \rightarrow y$ as $n \rightarrow \infty$ we have $x \in \mathcal{D}(F)$ and $F(x)=y$. 

\vspace*{2ex}

In the following, we want to revisit some basic definitions and statements for Banach spaces and their duals.

\subsection{Duality mappings and geometric properties of Banach spaces}
In \cite{ws16} we see that the (sub)gradient
\begin{displaymath}
 g(x) := \partial \Big(\frac{1}{2} \left\lVert F(\cdot) - y^{\delta} \right\rVert_Y^2 \Big)(x)
\end{displaymath}
of the least squares functional
\begin{displaymath}
 \Psi(x) := \frac{1}{2} \left\lVert F(x) - y^{\delta} \right\rVert_Y^2
\end{displaymath}
plays an essential role in the SESOP/RESESOP algorithms. This motivates the use of \emph{subgradients} of the more general convex functionals
\begin{displaymath}
 \Psi^p(x) := \frac{1}{p} \left\lVert F(x) - y^{\delta} \right\rVert_Y^p,
\end{displaymath}
where $p>1$, and we have
\begin{displaymath}
 \partial \Psi^p(x) = F'(x)^* \left( \Big( \partial \frac{1}{p} \left\lVert \cdot \right\rVert_Y^p \Big) \big(F(x) - y^{\delta}\big) \right).
\end{displaymath}
According to Asplund's Theorem (see, e.g., \cite{ea67}), we have 
\begin{equation}
 \partial \left( \frac{1}{p} \lVert \cdot \rVert_X^p \right) = J_p^X : X \rightrightarrows X^*,
\end{equation}
where $\big(X,\lVert \cdot \rVert_X\big)$ is a Banach space, $p \geq 1$, and
\begin{equation} \label{def_dualmap}
 J_p^X (x) := \left\lbrace \xs \in X^* \, : \, \langle \xs,x \rangle = \lVert \xs \rVert \cdot \lVert x \rVert, \, \lVert \xs \rVert = \lVert x \rVert^{p-1} \right\rbrace
\end{equation}
is the (set-valued) duality mapping of $X$ with gauge function $t \mapsto t^{p-1}$.

\vspace*{2ex}

\begin{remark} The following statements can be found in \cite{cioranescu12, skhk12}.
 \begin{itemize}
  \item[(a)] We obtain
  \begin{displaymath}
   \big(\partial \Psi^p\big) (x) = F'(x)^* J^Y_p \big(F(x) - y^{\delta}\big)
  \end{displaymath}
  for the subgradient of $\Psi^p$ in $x \in X$. 
  \item[(b)] If $X$ is a Hilbert space and $p=2$, we have $J_2^X(x)=x$.
  \item[(c)] The duality mappings are monotone, i.e., we have
  \begin{displaymath}
   \left\langle x^* - \widetilde{x}^* , x - \widetilde{x} \right\rangle \geq 0
  \end{displaymath}
   for all $x, \widetilde{x} \in X$, $x^* \in J^X_p(x)$ and $\widetilde{x}^* \in J^X_p(\widetilde{x})$.
 \end{itemize}
\end{remark}

\vspace*{2ex}

In order to obtain a better understanding of duality mappings, we have to take a closer look at the geometrical properties of Banach spaces, which influence the properties of the duality mappings. The notion of convexity and smoothness are especially important.

\vspace*{2ex}

\begin{definition}
 Let $X$ be a Banach space. We call $X$
 \begin{itemize}
  \item[(i)] \emph{strictly convex}, if 
   \begin{displaymath}
    \left\lVert \frac{1}{2}(x + \widetilde{x}) \right\rVert < 1 
   \end{displaymath}
   for all $x, \widetilde{x} \in \left\lbrace v \in X \, : \, \lVert v \rVert = 1 \right\rbrace$ with $x \neq \widetilde{x}$,
  \item[(ii)] \emph{smooth}, if for every $0 \neq x \in X$ there is a unique $\xs \in \Xs$ satisfying $\lVert \xs \rVert = 1$ and $\xs(x) = \lVert x \rVert$. 
 \end{itemize}

\end{definition}

\vspace*{2ex}

\begin{definition}
 Let $X$ be a Banach space. $X$ is \emph{$p$-convex}, if there is a constant $c_p>0$, such that
 \begin{equation} \label{p_convex}
  \frac{1}{p} \lVert x-\widetilde{x} \rVert^p \geq \frac{1}{p} \lVert x \rVert^p - \big\langle j_p^X(x), \widetilde{x} \big\rangle + \frac{c_p}{p} \lVert \widetilde{x} \rVert^p
 \end{equation}
 for all $x,\widetilde{x} \in X$ and all $j_p^X \in J^X_p$. \\
 We call $X$ \emph{$p$-smooth}, if there exists a constant $G_p>0$, such that
 \begin{equation} \label{p_smooth}
 \frac{1}{p} \lVert x-\widetilde{x} \rVert^p \leq \frac{1}{p} \lVert x \rVert^p - \big\langle j_p^X(x), \widetilde{x} \big\rangle + \frac{G_p}{p} \lVert \widetilde{x} \rVert^p
 \end{equation}
 for all $x,\widetilde{x} \in X$ and all $j_p^X \in J^X_p$. 
\end{definition}

\vspace*{2ex}

\begin{definition}
 Let $X$ be a Banach space. We call $X$ \emph{uniformly convex}, if the \emph{modulus of convexity}
 \begin{displaymath}
  \delta_X\, : \, [0,2]\rightarrow [0,1], \ \delta_X(\varepsilon) := \inf \left\lbrace 1-\left\lVert \frac{1}{2}(x+\widetilde{x}) \right\rVert \, : \, \lVert x \rVert = \lVert \widetilde{x} \rVert = 1, \, \lVert x-\widetilde{x} \rVert \geq \varepsilon \right\rbrace
 \end{displaymath}
 fulfills $\delta_X(\varepsilon) > 0$ for all $0 < \varepsilon \leq 2$.
\end{definition}

\vspace*{2ex}

\begin{definition}
 Let $X$ be a Banach space. We call $X$ \emph{uniformly smooth}, if the \emph{modulus of smoothness}
 \begin{displaymath}
  \rho_X \, : \, [0,\infty) \rightarrow [0,\infty), \ \rho_X(\tau) := \frac{1}{2}\sup\left\lbrace \lVert x+\widetilde{x} \rVert + \lVert x-\widetilde{x} \rVert - 2 \, : \, \lVert x \rVert = 1, \lVert \widetilde{x} \rVert \leq \tau \right\rbrace
 \end{displaymath}
 fulfills $\lim_{\tau \rightarrow 0} \tau^{-1} \cdot \rho_X(\tau) =0$.
\end{definition}

\vspace*{2ex}

\begin{remark}
 If the Banach space $X$ is $p$-smooth, then $\lVert \cdot \rVert_X^p$ is Fr\'echet differentiable and thus $J_p^X(x)$ is single-valued for all $x\in X$ (see, e.g., \cite{skhk12}). 
\end{remark}

\vspace*{2ex}

\begin{theorem} \label{theorem_properties_dm}
 The following statements hold.
\begin{itemize}
 \item[(i)] Let $X$ be a uniformly convex or uniformly smooth Banach space. Then $X$ is reflexive.
 \item[(ii)] If $X$ is a uniformly smooth Banach space, then $J_p^X$ is single-valued and uniformly continuous on bounded sets.
 \item[(iii)] If the Banach space $X$ is $s$-convex for some $s >0$ and smooth, then $J_p^X$ is single-valued, norm-to-weak continuous, bijective, and the duality mapping $J_{p^*}^{X^*}$ is single-valued. We then have
 \begin{equation}
  J_{p^*}^{X^*}\big( J_{p}^{X}(x) \big) = x.
 \end{equation}
\end{itemize}
\end{theorem}

\vspace*{2ex}

\begin{proof}
 See, e.g., \cite{sls06, skhk12} and the literature cited therein.
\end{proof}

\vspace*{2ex}

\begin{remark}
 If a Banach space $X$ is $p$-convex with $p\geq2$ ($p$-smooth with $p\leq 2$), then $X$ is also uniformly convex (uniformly smooth). In addition, uniform convexity (uniform smoothness) implies strict convexity (smoothness). Thus, $p$-convexity and $p$-smoothness are the "strongest" properties, whereas strict convexity and smoothness are the "weakest" properties (see, e.g., \cite{skhk12}). \\ 
 There is a range of further statements concerning the geometry of Banach spaces that we want to skip at this point. It is however noteworthy, that reflexive Banach spaces are an important class of Banach spaces with many useful geometric properties. In particular, the Lebesgue spaces $L^q(\Omega)$, $1<q<\infty$, which are often used in inverse problems are reflexive.
\end{remark}

\subsection{Bregman distances and Bregman projections}
In view of the reconstruction algorithms, it is convenient to employ \emph{Bregman distances} instead of the usual norm distances. The Bregman distance is later used to define the \emph{Bregman projection}. Our definition of the Bregman distance coincides with the definition in \cite{sls06, ss09, ssl08}.

\vspace*{2ex}

\begin{definition}
 Let $\Lambda\, : \, X \rightarrow \RR \cup \lbrace \infty \rbrace$ be a convex functional and $x,\widetilde{x} \in X$, $\lambda \in \partial\Lambda(x)$. The \emph{Bregman distance} of $x$ and $\widetilde{x}$ in $x$ and $\lambda$ w.r.t.~$\Lambda$ is defined as
 \begin{displaymath}
  D_{\lambda}^{\Lambda}(x,\widetilde{x}) := \Lambda(\widetilde{x}) - \Lambda(x) - \langle \lambda, \widetilde{x} - x \rangle_{\Xs \times X}. 
 \end{displaymath}
 Now let $1<p<\infty$ and $\Lambda := \Psi_p := \frac{1}{p} \lVert \cdot \rVert_X^p$, such that $\partial \Psi_p = J^X_p$. In this case, we write $D_p$ instead of $D_{j^X_p}$ for the Bregman distance, where
 \begin{equation} \label{breg1}
  D_{p}(x,\widetilde{x}) = \frac{1}{p} \lVert \widetilde{x} \rVert^p - \frac{1}{p} \lVert x \rVert^p - \left\langle j^X_p(x), \widetilde{x} - x\right\rangle
 \end{equation}
 for $x,\widetilde{x} \in X$, $j^X_p(x) \in J^X_p(x)$.
\end{definition}

\vspace*{2ex}

If $X$ is smooth, the duality mapping $J^X_p$ is single-valued according to Theorem \ref{theorem_properties_dm} (iii) and the Bregman distance $D_p(x,\widetilde{x})$ of $x,\widetilde{x} \in X$ can be calculated according to
\begin{align}
 D_{p}(x,\widetilde{x}) &= \frac{1}{p} \lVert \widetilde{x} \rVert^p + \frac{1}{p^*} \lVert x \rVert^p - \left\langle J^X_p(x), \widetilde{x} \right\rangle \label{breg2} \\
                        &= \frac{1}{p^*} \big( \lVert x \rVert^p - \lVert \widetilde{x} \rVert^p \big) + \left\langle J^X_p(\widetilde{x}) - J^X_p(x), \widetilde{x} \right\rangle. \label{breg3}
\end{align}

\vspace*{2ex}

\begin{theorem} \label{theorem_properties_bd}
 Let $X$ be a Banach space. For a fixed single-valued $j^X_p \in J^X_p$ and $x,\widetilde{x}\in X$ the corresponding Bregman distance has the following properties.
 \begin{itemize}
  \item[(i)] $D_{j^X_p}(x,\widetilde{x}) \geq 0$.
  \item[(ii)] $D_{j^X_p}(x,\widetilde{x}) = 0$ if and only if $j^X_p(\widetilde{x}) \in J^X_p(x)$.
  \item[(iii)] Let $X$ be smooth and uniformly convex. A sequence $\lbrace x_n \rbrace_{n\in\NN}$ in $X$ is bounded in $X$ if the sequence $\lbrace D_{j^X_p}(x_n,x) \rbrace$ is bounded in $\RR$.
  \item[(iv)] $D_{j^X_p}$ is continuous in its second argument. If $X$ is smooth and uniformly convex, $D_{j^X_p}$ is continuous also in its first argument.
  \item[(v)] If $X$ is smooth and uniformly convex, the statements
   \begin{itemize}
    \item[(a)] $\lim_{n\rightarrow\infty} \lVert x_n - x \rVert = 0$,
    \item[(b)] $\lim_{n\rightarrow\infty} \lVert x_n \rVert = \lVert x \rVert$ and $\lim_{n\rightarrow\infty} \langle j^X_p(x_n), x \rangle = \langle j^X_p(x), x \rangle$,
    \item[(c)] $\lim_{n\rightarrow\infty} D_{j^X_p}(x_n,x) = 0$
   \end{itemize}
   are equivalent.
  \item[(vi)] A sequence $\lbrace x_n \rbrace_{n\in\NN}$ in $X$ is a Cauchy sequence if it is bounded and for all $\varepsilon > 0$ there is an $N(\varepsilon) \in \NN$ such that $D_{j^X_p}(x_k,x_l) < \varepsilon$ for all $k,l \geq N(\varepsilon)$.
  \item[(vii)] $X$ is
   \begin{itemize}
    \item[(a)] $p$-convex if and only if $D_{j^X_p}(x,\widetilde{x}) \geq C \lVert x-\widetilde{x} \rVert^p$,
    \item[(b)] $p$-smooth if and only if $D_{j^X_p}(x,\widetilde{x}) \leq C \lVert x-\widetilde{x} \rVert^p$.
   \end{itemize}
 \end{itemize}
\end{theorem}

\vspace*{2ex}

\begin{proof}
 See \cite{skhk12, sls06, zxgr91}.
\end{proof}

\vspace*{2ex}

\begin{example} \label{ex1}
 Let $\Omega \subseteq \RR^N$, $N \in \NN$, denote an open domain. The Lebesgue spaces $L^q(\Omega)$ and the Sobolev spaces $W^{m,q}(\Omega)$, equipped with the respective usual norm, are $\max\lbrace 2,q \rbrace$-convex and $\min\lbrace 2,q \rbrace$-smooth for $1 < q < \infty$. In particular, they are  uniformly convex and uniformly smooth, such that the duality mappings $J_p$ are single-valued and continuous.
\end{example}

\vspace*{2ex}

\begin{definition}
 Let $\emptyset \neq C \subseteq X$ be a closed, convex set. The \emph{Bregman projection} of $x \in X$ onto $C$ with respect to $\Psi_p$ is the unique element $\Pi_C^p(x) \in C$, such that
 \begin{equation}
  D_p(x,\Pi_C^p(x)) = \min_{z \in C} D_p(x,z).
 \end{equation}
\end{definition}

\vspace*{2ex}

\begin{remark}
 Alternatively, we can characterize the Bregman projection of an element $x \in X$ onto $C$ w.r.t.~$\Psi_p$ as the element $\widetilde{x} \in C$, which fulfills the variational inequality
 \begin{displaymath}
  \left\langle J^X_p(\widetilde{x}) - J^X_p(x), z-\widetilde{x}  \right\rangle \geq 0 \quad \text{for all} \ z \in C.
 \end{displaymath}
 The latter is equivalent to the \emph{descent property}
 \begin{equation} \label{desc_prop}
  D_p(\widetilde{x},z) \leq D_p(x,z) - D_p(x,\widetilde{x}) \quad \text{for all} \ z \in C. 
 \end{equation}
\end{remark}

\vspace*{2ex}

\begin{definition}
 The open ball centered about $x_0 \in X$ with radius $\varrho > 0$ with respect to the Bregman distance $D_{j^X_p}$ is defined as
 \begin{displaymath}
  B_{\varrho}^{p}(\widetilde{x}) := B_{\varrho}^{D_{j^X_p}}(\widetilde{x}) := \left\lbrace x \in X \, : \, D_{j^X_p}(x,\widetilde{x}) \leq \varrho \right\rbrace.
 \end{displaymath}
\end{definition}

\vspace*{2ex}

Note that this definition is not symmetric, i.e., the center of the open ball is the second argument in the Bregman distance.

\subsection{Hyperplanes, halfspaces, and stripes}
The overall idea is to describe an iterative method where we approximate a solution of an inverse problem by sequentially projecting the iterates onto subsets (hyperplanes, stripes) of the source space $X$ that contain the solution set. 

\vspace*{2ex}

\begin{definition}
 Let $0 \neq u^* \in X^*$ and $\alpha, \xi \in \RR$ with $\xi \geq 0$. We define the \emph{hyperplane}
 \begin{displaymath}
  H(u^*,\alpha) := \left\lbrace x\in X \, : \, \left\langle u^*,x \right\rangle = \alpha \right\rbrace
 \end{displaymath}
 and the \emph{halfspace}
 \begin{displaymath}
  H_{\leq}(u^*,\alpha) := \left\lbrace x\in X \, : \, \left\langle u^*,x \right\rangle \leq \alpha \right\rbrace.
 \end{displaymath}
 Analogously, we define $H_{\geq}(u^*,\alpha)$, $H_{<}(u^*,\alpha)$ and $H_{>}(u^*,\alpha)$. Finally, the set
 \begin{displaymath}
  H(u^*,\alpha, \xi) := \left\lbrace x\in X \, : \, \left\lvert \left\langle u^*,x \right\rangle - \alpha \right\rvert \leq \xi \right\rbrace.
 \end{displaymath}
 is called a \emph{stripe}.
\end{definition}

\vspace*{2ex}

Note that the sets $H(u^*,\alpha)$, $H_{\leq}(u^*,\alpha)$, $H_{\geq}(u^*,\alpha)$, and $H(u^*,\alpha, \xi)$ are closed, convex, nonempty subsets of $X$, whereas $H_{<}(u^*,\alpha)$ and $H_{>}(u^*,\alpha)$ are not closed.

\subsection{Bregman projections onto hyperplanes}
The Bregman projection of $x\in X$ onto a hyperplane or intersections of hyperplanes can be described as a convex optimization problem. The following results are essential for our techniques. All statements have been shown in the articles \cite{ss09, ssl08} by Sch\"opfer, Schuster, and Louis for linear operator equations.

\vspace*{2ex}

\begin{proposition} \label{prop_prop_h}
 Let $X$ be reflexive, smooth and uniformly convex, such that $J_p^X$ and $J_{p^*}^{X^*}$ are single-valued. Let $H(u_1^*, \alpha_1),...,H(u_N^*, \alpha_N)$ denote $N\in \NN$ hyperplanes with nonempty intersection
 \begin{displaymath}
  H:=\bigcap_{k=1}^N H(u_k^*, \alpha_k).
 \end{displaymath}
 The Bregman projection $\Pi_H^p(x)$ of $x\in X$ onto $H$ is given by
 \begin{equation}
  \Pi_H^p(x) = J^{\Xs}_{p^*} \left( J^X_p(x) - \sum_{k=1}^N \widetilde{t}_k u_k^* \right),
 \end{equation}
 where the vector $\widetilde{t} = (\widetilde{t}_k)_{k=1,...,N} \in \RR^N$ solves the minimization problem
 \begin{equation}
  \min_{t \in \RR^N} h(t) := \frac{1}{p^*} \left\lVert J^X_p(x) - \sum_{k=1}^N t_k u_k^* \right\rVert^{p^*} + \sum_{k=1}^N t_k \alpha_k.
 \end{equation}
 The function $h:\RR^N \rightarrow \RR$ is convex and has continuous partial derivatives
 \begin{displaymath}
  \partial_{t_j} h(t) = \left\langle u_j^*, J^{\Xs}_{p^*} \left( J^X_p(x) - \sum_{k=1}^N t_k u_k^* \right) \right\rangle + \alpha_j, \quad j=1,...,N.
 \end{displaymath}
 If $u_1^*,...,u_N^*$ are linearly independent, $h$ is strictly convex and $\widetilde{t}$ is unique.
\end{proposition}

\vspace*{2ex}

\begin{proposition} \label{prop_proj}
 The following statements are helpful for the analysis of our methods.
 \begin{itemize}
  \item[(i)] Consider two halfspaces $H_1 := H_{\leq}(u_1^*,\alpha_1)$ and $H_2 := H_{\leq}(u_2^*,\alpha_2)$ with linearly independent vectors $u_1^*, u_2^*$. Then we have
   \begin{displaymath}
    \widetilde{x} = \Pi^p_{H_1 \cap H_2}(x)
   \end{displaymath}
   if and only if $\widetilde{x}$ fulfills the Karush-Kuhn-Tucker (KKT) conditions for the minimization problem
   \begin{displaymath}
    \min_{z\in H_1\cap H_2} D_p(x,z),
   \end{displaymath}
   which read
   \begin{equation}
    \begin{split}
     J^X_p(\widetilde{x}) = J^X_p(x) - t_1 u_1^* - t_2 u_2^*, \quad t_1, t_2 \geq 0, \\
     \langle u_1^*, \widetilde{x} \rangle \leq \alpha_1, \ \ \langle u_2^*, \widetilde{x} \rangle \leq \alpha_2, \\
     t_1\big( \alpha_1 - \langle u_1^*, \widetilde{x} \rangle \big) \leq 0, \ \ t_2\big( \alpha_2 - \langle u_2^*, \widetilde{x} \rangle \big) \leq 0.
     \end{split}
   \end{equation}
  \item[(ii)] For $x \in H_>(u^*,\alpha)$ the Bregman projection of $x$ onto $H_{\leq}(u^*,\alpha)$ is given by
   \begin{equation}
    \Pi^p_{H_{\leq}(u^*,\alpha)}(x) = \Pi^p_{H(u^*,\alpha)}(x) = J^{\Xs}_{p^*} \big( J^X_p(x) - t_+u^* \big),
   \end{equation}
   where $t_+ > 0$ is the unique, positive solution of
   \begin{equation}
    \min_{t\in\RR} \frac{1}{p^*} \left\lVert J^X_p(x) - tu^* \right\rVert^{p^*} + \alpha t.
   \end{equation}
  \item[(iii)] The Bregman projection of $x \in X$ onto a stripe $H(u^*,\alpha,\xi)$ is given by
   \begin{equation}
    \Pi^p_{H(u^*,\alpha,\xi)}(x) = \begin{cases}
                                    \Pi^p_{H_{\leq}(u^*,\alpha+\xi)}(x), & x \in H_{>}(u^*,\alpha+\xi), \\ 
                                    x,                                   & x \in H(u^*,\alpha,\xi), \\ 
                                    \Pi^p_{H_{\geq}(u^*,\alpha-\xi)}(x), & x \in H_{<}(u^*,\alpha-\xi).
                                   \end{cases}
   \end{equation}
 \end{itemize}
\end{proposition}

\section{Sequential subspace optimization in Banach spaces} \label{sec_sesop}
We have assembled all the tools that are necessary to formulate the SESOP and RESESOP algorithm for nonlinear inverse problems in Banach spaces. The first algorithm represents a general procedure for exact data, whereas Algorithm \ref{algo2} is designed for noisy data. Finally, we present a special case of Algorithm \ref{algo2} by specifying the search directions. 

\subsection{SESOP for exact data}
We begin with some definitions, before we introduce the general SESOP method for exact data.

\vspace*{2ex}

\begin{definition} \label{def_param_stripes_1}
 For $n\in\NN$ and a finite index set $I_n\subseteq \lbrace 0,1,...,n \rbrace$ we define the stripes
 \begin{displaymath}
  H_{n,i} := H\big(u^*_{n,i}, \alpha_{n,i}, \xi_{n,i}\big)
 \end{displaymath}
 for all $i\in I_n$, where
 \begin{displaymath}
 \begin{split}
   u^*_{n,i} &:= F'(x_i)^*w_{n,i}, \\
   \alpha_{n,i} &:= \big\langle F'(x_i)^*w_{n,i}, x_i \big\rangle_{X^* \times X} - \big\langle w_{n,i}, R_i \big\rangle_{Y^*\times Y} \\
   \xi_{n,i} &:= \ctc \lVert w_{n,i} \rVert \cdot \lVert R_i \rVert.
 \end{split}
 \end{displaymath}
 We denote the residual in $x_i$ by
 \begin{displaymath}
  R_i := F(x_i) - y.
 \end{displaymath}
\end{definition}

\vspace*{2ex}

\begin{algorithm} \label{algo1}
 Choose an initial value $x_0 \in X$. At iteration $n \in \NN$, choose a finite index set $I_n$ and search directions $u^*_{n,i} := F'(x_i)^*w_{n,i}$ with $w_{n,i} \in Y^*$. For each $i \in I_n$, we define the stripes
  \begin{displaymath}
  H_{n,i} := H\big(u^*_{n,i}, \alpha_{n,i}, \xi_{n,i}\big)
 \end{displaymath}
 as in Definition \ref{def_param_stripes_1} and compute the new iterate as
 \begin{displaymath}
  x_{n+1} = J_{p^*}^{X^*} \left(J^X_p(x_n) - \sum_{i\in I_n} t_{n,i} F'(x_i)^*w_{n,i} \right),
 \end{displaymath}
 where $t_n := (t_{n,i})_{i\in I_n}$ is chosen such that
 \begin{displaymath}
  x_{n+1} = \Pi_{H_n}^p (x_n), \quad H_n := \bigcap_{i \in I_n} H_{n,i}
 \end{displaymath}
 and
 \begin{displaymath}
  D_p(x_{n+1},z) \leq D_p(x_n,z) - C\cdot \lVert R_n \rVert^p
 \end{displaymath}
 for all $z \in M_{F(x)=y}$, where $C = C(c_F, \ctc) > 0$.
\end{algorithm}

\vspace*{2ex}

\begin{remark}
 In Algorithm \ref{algo1}, the new iterate $x_{n+1}$ is calculated as the Bregman projection of the current iterate $x_n$ onto the intersection of the stripes $H_{n,i}$. At this point, the optimality conditions from Proposition \ref{prop_proj} can be used to determine the parameters $t_{n,i}$.
\end{remark}

\subsection{RESESOP for noisy data}
We now adapt Algorithm \ref{algo1} to the case of noisy data $y^{\delta}$. To this end, we mainly have to increase the width of the stripes in accordance with the noise level $\delta$. 

\vspace*{2ex}

\begin{definition} \label{def_param_stripes_2}
 For $n\in\NN$ and a finite index set $I_n^{\delta}\subseteq \lbrace 0,1,...,n \rbrace$ we define the stripes
 \begin{displaymath}
  H_{n,i}^{\delta} := H\big(u^*_{n,i}, \alpha^{\delta}_{n,i}, \xi^{\delta}_{n,i}\big)
 \end{displaymath}
 for all $i\in I_n^{\delta}$, where
 \begin{displaymath}
 \begin{split}
   u^*_{n,i} &:= F'(x_i^{\delta})^*w^{\delta}_{n,i}, \\
   \alpha^{\delta}_{n,i} &:= \big\langle F'(x^{\delta}_i)^*w^{\delta}_{n,i}, x^{\delta}_i \big\rangle_{X^* \times X} - \big\langle w^{\delta}_{n,i}, R^{\delta}_i \big\rangle_{Y^*\times Y} \\
   \xi^{\delta}_{n,i} &:= \left(\delta + \ctc \left( \big\lVert R^{\delta}_i \big\rVert + \delta \right) \right) \lVert w^{\delta}_{n,i} \rVert.
 \end{split}
 \end{displaymath}
 As before, we denote the residual in $x_i^{\delta}$ by
 \begin{displaymath}
  R^{\delta}_i := F\big(x^{\delta}_i\big) - y^{\delta}.
 \end{displaymath}
\end{definition}

\vspace*{2ex}

\begin{algorithm} \label{algo2}
 Choose an initial value $x^{\delta}_0 := x_0\in X$. At iteration $n \in \NN$, choose a finite index set $I^{\delta}_n$ and search directions $u^*_{n,i} := F'\big(x^{\delta}_i\big)^*w^{\delta}_{n,i}$ with $w^{\delta}_{n,i} \in Y^*$. For each $i \in I^{\delta}_n$, we define the stripes
  \begin{displaymath}
  H^{\delta}_{n,i} := H\big(u^*_{n,i}, \alpha^{\delta}_{n,i}, \xi^{\delta}_{n,i}\big)
 \end{displaymath}
 as in Definition \ref{def_param_stripes_1} and compute the new iterate as
 \begin{displaymath}
  x^{\delta}_{n+1} = J_{p^*}^{X^*} \left(J^X_p\big(x^{\delta}_n\big) - \sum_{i\in I^{\delta}_n} t^{\delta}_{n,i} F'\big(x^{\delta}_i\big)^*w^{\delta}_{n,i} \right),
 \end{displaymath}
 where $t^{\delta}_n := \big(t^{\delta}_{n,i}\big)_{i\in I^{\delta}_n}$ are chosen such that
 \begin{displaymath}
  x^{\delta}_{n+1} = \Pi_{H^{\delta}_n}^p \big(x^{\delta}_n\big) \in H^{\delta}_n := \bigcap_{i \in I^{\delta}_n} H^{\delta}_{n,i}
 \end{displaymath}
 and
 \begin{displaymath}
  D_p\big(x^{\delta}_{n+1},z\big) \leq D_p\big(x^{\delta}_n,z\big) - C\cdot \big\lVert R^{\delta}_n \big\rVert^p
 \end{displaymath}
 for all $z \in M_{F(x)=y}$, where $C = C(c_F, \ctc) > 0$.
\end{algorithm}

\vspace*{2ex}

\begin{proposition}
 For all $n \in \NN$ and $i \in I_n$, resp.~$i \in I_n^{\delta}$, the stripes $H_{n,i}$ and $H_{n,i}^{\delta}$ contain the solution set $M_{F(x)=y}$.
\end{proposition}

\vspace*{1ex}

\begin{proof}
 We prove the statement for the case of noisy data, since it includes the special case $\delta=0$. For $z\in M_{F(x)=y}$ we have
 \begin{displaymath}
  \begin{split}
     \left\lvert \left\langle u^*_{n,i}, z \right\rangle - \alpha^{\delta}_{n,i} \right\rvert &= \left\lvert \left\langle F'(x_i^{\delta})^* w^{\delta}_{n,i}, z \right\rangle - \big(\big\langle F'(x^{\delta}_i)^*w^{\delta}_{n,i}, x^{\delta}_i \big\rangle - \big\langle w^{\delta}_{n,i}, R^{\delta}_i \big\rangle \big) \right\rvert \\
     &= \left\lvert \left\langle F'(x_i^{\delta})^*w^{\delta}_{n,i}, z - x^{\delta}_i \right\rangle + \big\langle w^{\delta}_{n,i}, R^{\delta}_i \big\rangle \right\rvert \\
     &= \left\lvert \left\langle w^{\delta}_{n,i}, F'(x_i^{\delta}) (z - x^{\delta}_i) + R^{\delta}_i \right\rangle \right\rvert \\
     &\leq \left\lVert w^{\delta}_{n,i} \right\rVert \cdot \left\lVert F(x_i^{\delta}) - F(z) - F'(x_i^{\delta}) (x^{\delta}_i - z) + y - y^{\delta}\right\rVert \\
     &\leq \left\lVert w^{\delta}_{n,i} \right\rVert \cdot \big( \ctc \left\lVert F(x_i^{\delta}) - F(z) \right\rVert + \delta \big) \\
     &\leq \left\lVert w^{\delta}_{n,i} \right\rVert \cdot \big( \ctc \big(\left\lVert F(x_i^{\delta}) - y^{\delta} \right\rVert + \delta \big) + \delta \big) = \xi^{\delta}_{n,i},
  \end{split}
 \end{displaymath}
 using the definitions of the parameters, the Cauchy-Schwarz inequality, and the tangential cone condition. Hence, $z \in H_{n,i}^{\delta}$.
\end{proof}

\section{Convergence and regularization} \label{sec_convreg}
This section is dedicated to an analysis of the methods we proposed in Section \ref{sec_sesop}. Previous results from \cite{ss09, ssl08, ws16} show that the use of the \emph{current gradient}
\begin{equation} \label{current_gradient}
 g_n^{\delta} := \partial \left( \frac{1}{2} \big\lVert F(\cdot) - y^{\delta} \big\rVert^2\right) (x^{\delta}_n) = F'\big(x_n^{\delta}\big)^* J^Y_2 \big( F\big(x_n^{\delta}\big) - y^{\delta} \big)
\end{equation}
assures the validity of the descent properties. For this choice, we present convergence results. By specifying the search directions even further, we are also able to prove that the RESESOP method yields a regularization technique.
Throughout this section we assume that the Banach space $X$ is uniformly smooth and $p$-convex. Consequently, the dual $\Xs$ of $X$ is $p^*$-smooth and uniformly convex, and the duality mappings $J^X_p: X \rightarrow X^*$, $J^{\Xs}_{p^*}: X^* \rightarrow X$ are single-valued and continuous. The data space $Y$ is an arbitrary Banach space. Additional properties are only required in Section \ref{subsec_resesop_reg}.

\subsection{Convergence of the SESOP algorithm}
We state some results for the SESOP algorithm with the following specifications, which assure that the current gradient is always contained in the search space. 

\vspace*{2ex}

\begin{assumption} \label{spec_1}
 \begin{itemize}
  \item[(i)] Let $n \in I_n$ for all iterations $n \in \NN$.
  \item[(ii)] Let $j_2^Y \in J_2^Y$. Set $w_{n,n} := j^Y_2 \big( F\big(x_n\big) - y \big)$ for all $n\in \NN$.
  \item[(iii)] Choose $I_n \subseteq \lbrace n-N+1,...,n \rbrace$ for some fixed $N>1$ such that the set of search directions
   \begin{displaymath}
    U^*_n :=\left\lbrace u^*_{n,i} := F'(x_i)^*w_{n,i} \, : \, i \in I_n \right\rbrace
   \end{displaymath}
   is linearly independent.
  \item[(iv)] The initial value $x_0$ fulfills $x_0 \in B^{p}_{\varrho}(z)$ for some $z \in M_{F(x)=y}$ and $\varrho > 0$.
 \end{itemize}
\end{assumption}

\vspace*{2ex}

\begin{proposition} \label{prop_descent_prop_algo1}
 Let $\lbrace x_n \rbrace_{n\in\NN}$ denote the sequence of iterates that is generated by Algorithm \ref{algo1} with the specifications \ref{spec_1}. We then have
 \begin{equation}\label{xn_above_Hnn}
  x_n \in H_>\big(u^*_{n,n},\alpha_{n,n} + \xi_{n,n}\big).
 \end{equation}
 By projecting $x_n$ first onto $H_{n,n} := H\big(u^*_{n,n},\alpha_{n,n}, \xi_{n,n}\big)$, we obtain the descent property
 \begin{equation} \label{descent_1}
  D_p(x_{n+1},z) \leq D_p\big(\Pi^p_{H_{n,n}} \big(x_{n}\big),z \big) \leq D_p\big( x_n,z \big) - \frac{(1-\ctc)^p}{p \cdot G_{p^*}^{p-1} \cdot c_F^p} \big\lVert R_n \big\rVert^p 
 \end{equation}
 for $z \in M_{F(x)=y}$, which implies $x_n \in B^{p}_{\varrho}(z)$.
\end{proposition}

\vspace*{2ex}

\begin{proof}
 To prove \eqref{xn_above_Hnn} we set $w_{n,n} := j_2^Y \big( F\big(x_n\big) - y \big)$ and use
 \begin{displaymath}
  \left\langle j_2^Y \big( F\big(x_n\big) - y \big), F\big(x_n\big) - y \right\rangle = \big\lVert F\big(x_n\big) - y \big\rVert^2
 \end{displaymath}
 according to the definition of $j_2^Y$, see \eqref{def_dualmap}. We estimate
 \begin{displaymath}
  \alpha_{n,n} + \xi_{n,n} = \left\langle u_{n,n}^*, x_n \right\rangle - (1-\ctc) \left\lVert F(x_n) - y \right\rVert^2 \leq \left\langle u_{n,n}^*, x_n \right\rangle
 \end{displaymath}
 due to $0 < \ctc < 1$, which yields \eqref{xn_above_Hnn}. \\
 In order to obtain the descent property \eqref{descent_1} for the Bregman distance, we exploit Proposition \ref{prop_proj}. First of all, we note that due to \eqref{xn_above_Hnn}, we have
 \begin{displaymath}
  \Pi^p_{H(u_{n,n}^*,\alpha_{n,n},\xi_{n,n})}(x_n) = \Pi^p_{H(u_{n,n}^*,\alpha_{n,n}+\xi_{n,n})}(x_n).
 \end{displaymath}
 For $z \in M_{F(x)=y}$ we thus obtain 
 \begin{displaymath}
 \begin{split}
  D_p\big(x_{n+1},z\big) &\leq D_p\big(\Pi^p_{H_{n,n}}(x_n),z\big) \\
                         &= \frac{1}{p^*} \big\lVert J^X_p(x_n) - t_+ u_{n,n}^* \big\rVert^{p^*} + t_+ \left\langle u_{n,n}^*, z \right\rangle - \left\langle J^X_p(x_n),z \right\rangle + \frac{1}{p} \lVert z \rVert^p \\
                         &\leq \frac{1}{p^*} \big\lVert J^X_p(x_n) - t_+ u_{n,n}^* \big\rVert^{p^*} + t_+ \big( \alpha_{n,n} + \xi_{n,n} \big) - \left\langle J^X_p(x_n),z \right\rangle + \frac{1}{p} \lVert z \rVert^p \\
                         &= h(t_+) - \left\langle J^X_p(x_n),z \right\rangle + \frac{1}{p} \lVert z \rVert^p \\
                         &\leq h(\widetilde{t}) - \left\langle J^X_p(x_n),z \right\rangle + \frac{1}{p} \lVert z \rVert^p,
 \end{split}
 \end{displaymath}
 where $t_+$ minimizes the function
 \begin{displaymath}
  h(t) = \frac{1}{p^*} \big\lVert J^X_p(x_n) - t\cdot u_{n,n}^* \big\rVert^{p^*} + t \cdot \big( \alpha_{n,n} + \xi_{n,n} \big),
 \end{displaymath}
 and for
 \begin{displaymath}
  \widetilde{t} := \left( \frac{\big\langle u_{n,n}^*, x_n \big\rangle - \big( \alpha_{n,n} + \xi_{n,n} \big)}{G_{p^*} \big\lVert u_{n,n}^* \big\rVert^{p^*}} \right)^{p-1} > 0
 \end{displaymath}
 we thus have $h(t_+) \leq h(\widetilde{t})$, which we have used in the last estimate. Note that the choice of $\widetilde{t}$ is inspired by the well-known explicit form of the optimization parameter $t_+$ in the Hilbert space case, see \cite{ws16} and also \cite{skhk12}. \\
 We now estimate the value of $h$ in $\widetilde{t}$. Since $X$ is $p$-convex, $X^*$ is $p^*$-smooth and \eqref{p_smooth} applies. Together with \eqref{def_dualmap} (and by plugging in $\widetilde{t}$) we obtain
 \begin{displaymath}
  \begin{split}
   h(\widetilde{t}) &\leq \frac{1}{p^*} \big\lVert J^X_p(x_n) \big\rVert^{p^*} - \left\langle J^{X^*}_{p^*}\big( J^X_p(x_n) \big), \widetilde{t} u_{n,n}^* \right\rangle + \frac{G_{p^*}}{p^*} \big\lVert \widetilde{t} u_{n,n}^* \big\rVert^{p^{*}} + \widetilde{t}\cdot \big( \alpha_{n,n} + \xi_{n,n} \big) \\
                    &= \frac{1}{p^*} \big\lVert x_n \big\rVert^{p} - \frac{G_{p^*}}{p} \left( \frac{\big\langle u_{n,n}^*, x_n \big\rangle - \big( \alpha_{n,n} + \xi_{n,n} \big)}{G_{p^*} \big\lVert u_{n,n}^* \big\rVert} \right)^p.
  \end{split}
 \end{displaymath}
 Inserting this in the previous estimate and using \eqref{breg2}, we arrive at
 \begin{displaymath}
  \begin{split}
  D_p\big(x_{n+1},z\big) & \leq D_p\big(x_n,z\big) - \frac{1}{p G_{p^*}^{p-1}} \left( \frac{\big\langle u_{n,n}^*, x_n \big\rangle - \big( \alpha_{n,n} + \xi_{n,n} \big)}{ \big\lVert u_{n,n}^* \big\rVert} \right)^p \\
                         & \leq D_p\big(x_n,z\big) - \frac{(1-\ctc)^p}{p G_{p^*}^{p-1}\cdot c_F^p} \cdot \lVert R_n \rVert^p.
  \end{split}
 \end{displaymath}
 In the last step we have used the definitions of $u_{n,n}^*$, $\alpha_{n,n}$ and $\xi_{n,n}$ as well as the properties of $F$ we postulated in Section \ref{sec_prelim}, together with $\big\lVert F'(x_n)^* \big\rVert = \big\lVert F'(x_n) \big\rVert \leq c_F$.
\end{proof}

\vspace*{2ex}

\begin{proposition}\label{prop_weak_conv_algo1}
 The sequence of iterates $\lbrace x_n \rbrace_{n\in\NN}$ that is generated by Algorithm \ref{algo1} under the assumptions \ref{spec_1} is bounded and has weak cluster points that solve \eqref{op_eq}. The residual fulfills
 \begin{displaymath}
  \lim_{n\rightarrow \infty} \lVert R_n \rVert = 0.
 \end{displaymath}
\end{proposition}

\vspace*{2ex}

\begin{proof}
 According to Proposition \ref{prop_descent_prop_algo1}, the sequence $\left\lbrace D_p(x_n,z) \right\rbrace_{n\in\NN}$ of Bregman distances is monotonically decreasing and bounded. As a consequence, the sequence of iterates $\lbrace x_n \rbrace_{n \in \NN}$ is also bounded, see Theorem 2.12 from \cite{sls06}, and has weak cluster points. From \eqref{descent_1} we obtain
 \begin{displaymath}
 \begin{split}
  \frac{(1-\ctc)^p}{p G_{p^*}^{p-1} c_{\mathrm{F}}^p}\sum_{n=0}^{M} \lVert R_{n} \rVert^p &\leq \sum_{n=0}^{M} \big( D_p(x_{n},z) - D_p(x_{n+1},z) \big) \\
                                                                                          &= D_p(x_{0},z) - D_p(x_{M+1},z) \\
                                                                                          &\leq D_p(x_{0},z)
 \end{split}
 \end{displaymath}
 for all $M\geq0$ and thus
 \begin{displaymath}
  \sum_{n=0}^{\infty} \lVert R_{n} \rVert^p \leq \frac{p G_{p^*}^{p-1}}{(1-\ctc)^p} \cdot D_p(x_{0},z),
 \end{displaymath}
 which implies $\lVert R_n \rVert \rightarrow 0$ for $n \rightarrow \infty$. \\
 Now let $\lbrace x_{n_k} \rbrace_{k \in \NN}$ denote a weakly convergent subsequence with $x_{n_k} \rightharpoonup: \hat{x}$ for $k \rightarrow \infty$. The sequence of residuals $\lVert R_{n_k} \rVert_{k\in\NN}$ is a null sequence. Due to the weak sequential closedness of $F$ and the continuity of the norm we have
 \begin{displaymath}
  0 =\lim_{k\rightarrow \infty} \lVert R_{n_k} \rVert = \lim_{k\rightarrow \infty} \lVert F(x_{n_k}) - y \rVert = \lVert F(\hat{x}) - y \rVert
 \end{displaymath}
 and consequently $\hat{x} \in M_{F(x)=y}$.
\end{proof}

\vspace*{2ex}

\begin{remark}
 In the situation of Proposition \ref{prop_weak_conv_algo1} we obtain a strongly converging subsequence, if $X$ is a finite dimensional space, since all weakly convergent sequences in a finite dimensional space are already strongly convergent.
\end{remark}

\vspace*{2ex}

\begin{theorem} \label{theorem_conv_exact_data}
 Let $w_{n,i}:=j_2^Y(R_i)$ for all $n\in\NN$ and $i\in I_n$ in addition to the specifications of Assumption \ref{spec_1}. The sequence of iterates $\lbrace x_n \rbrace_{n\in\NN}$ converges strongly to a solution $x^+$ of \eqref{op_eq}, if
 \begin{displaymath}
  \left\lvert t_{n,i} \right\rvert < t \quad \text{for all} \ n \in \NN \ \text{and} \ i \in I_n 
 \end{displaymath}
 for some $t>0$.
\end{theorem}

\vspace*{2ex}

\begin{proof}
 By the boundedness of $\lbrace x_n \rbrace_{n \in \NN}$, Proposition \ref{prop_weak_conv_algo1} yields the existence of a weakly convergent subsequence $\lbrace x_{n_k} \rbrace_{k \in \NN}$ with $x_{n_k} \rightharpoonup x^+$ for $k \rightarrow \infty$, such that $\left\lbrace \lVert x_{n_k} \rVert \right\rbrace$ converges. \\
 We now show that $\lbrace x_{n_k} \rbrace_{k \in \NN}$ is a Cauchy sequence. To this end, we use Theorem \ref{theorem_properties_bd} (vi) and consider the Bregman distance
 \begin{displaymath}
  D_p(x_{n_k},x_{n_l}) = \frac{1}{p^*} \left( \lVert x_{n_l} \rVert^p - \lVert x_{n_k} \rVert^p \right) + \left\langle J^X_p(x_{n_k}) - J^X_p(x_{n_l}), x_{n_k} \right\rangle
 \end{displaymath}
 of two iterates $x_{n_k}$ and $x_{n_l}$, where $k>l$. We have chosen the subsequence $\lbrace x_{n_k} \rbrace_{k \in \NN}$ such that the sequence $\lbrace\lVert x_{n_k} \rVert\rbrace_{k\in\NN}$ converges. We thus have
 \begin{displaymath}
  \frac{1}{p^*} \left( \lVert x_{n_l} \rVert^p - \lVert x_{n_k} \rVert^p \right) \rightarrow 0 \quad  \text{for} \ l \rightarrow \infty.
 \end{displaymath}
 The second term is reformulated such that
 \begin{displaymath}
 \begin{split}
  &\left\langle J^X_p(x_{n_k}) - J^X_p(x_{n_l}), x_{n_k} \right\rangle \\
  &\quad= \left\langle J^X_p(x_{n_k}) - J^X_p(x_{n_l}), x_{n_k} - x^+ \right\rangle + \left\langle J^X_p(x_{n_k}) - J^X_p(x_{n_l}), x^+ \right\rangle,
 \end{split}
 \end{displaymath}
 where the second term converges to 0 for $l \rightarrow \infty$, since the sequence $\lbrace x_{n_k} \rbrace_{k \in \NN}$ is weakly convergent due to Proposition \ref{prop_weak_conv_algo1}. The absolute value of the first term is estimated by using similar arguments as in the Hilbert space setting (\cite{ws16}, Theorem 4.4, and \cite{hns_cl}, Theorem 2.3). In particular, we make use of the recursion for the iterates $x_{n_j}$, $j=l+1,...,k$, and obtain
 \begin{displaymath}
  \begin{split}
   &\left\lvert\left\langle J^X_p(x_{n_k}) - J^X_p(x_{n_l}), x_{n_k} - x^+ \right\rangle\right\rvert \\
   &\quad= \left\lvert\left\langle - \sum_{j=n_l}^{n_k-1} \sum_{i\in I_j} t_{j,i} F'(x_j)^*j_2^Y(F(x_j)-y), x_{n_k} - x^+ \right\rangle\right\rvert \\
   &\quad\leq t \cdot \sum_{j=n_l}^{n_k-1} \sum_{i\in I_j} \left\lVert j_2^Y(F(x_j)-y) \right\rVert \cdot \left\lVert F'(x_j)\big(x_{n_k} - x^+\big) \right\rVert.
  \end{split}
 \end{displaymath}
 Without loss of generality we assume that the sequence $\left\lbrace \lVert F(x_{n_k}) - y \rVert \right\rbrace_{k\in\NN}$ is monotonically decreasing (otherwise, choose a monotonically decreasing subsequence), yielding
 \begin{displaymath}
  \left\lVert F(x_{n_k}) - y \right\rVert \leq \left\lVert F(x_j) - y \right\rVert =: \lVert R_j \rVert.
 \end{displaymath}
 Furthermore, we estimate
 \begin{displaymath}
  \begin{split}
   \left\lVert F'(x_j)\big(x_{n_k} - x^+\big) \right\rVert &= \left\lVert F'(x_j)\big(x_{n_k} - x_j + x_j - x^+\big) + F(x_i) - F(x_i) \right.\\
                                                           & \qquad \left. + F(x_{n_k}) - F(x_{n_k}) + y - y \right\rVert \\
                                                           &\leq \left\lVert F(x_j) - F(x_{n_k})- F'(x_j)\big(x_j - x_{n_k} \big) \right\rVert \\
                                                           & \qquad +\left\lVert F(x_j) - F(x^+)- F'(x_j)\big(x_j - x^+ \big) \right\rVert + \left\lVert F(x_{n_k}) - y \right\rVert \\
                                                           &\leq \ctc \left\lVert F(x_j) - F(x_{n_k}) \right\rVert + \ctc \left\lVert F(x_j) - F(x^+) \right\rVert + \left\lVert F(x_{n_k}) - y \right\rVert \\
                                                           &\leq 2\ctc\left\lVert F(x_j) - y \right\rVert + (1+\ctc) \left\lVert F(x_{n_k}) - y \right\rVert \\
                                                           &\leq (1+3\ctc) \lVert R_j \rVert
  \end{split}
 \end{displaymath}
 Together with the definition of the duality mapping $j_2^Y$, we arrive at
 \begin{displaymath}
  \left\lvert\left\langle J^X_p(x_{n_k}) - J^X_p(x_{n_l}), x_{n_k} - x^+ \right\rangle\right\rvert \leq t \cdot (1+3\ctc)\sum_{j=n_l}^{n_k-1} \sum_{i\in I_j} \lVert R_j \rVert^2.
 \end{displaymath}
 The right-hand side of the above estimate converges to zero for $l \rightarrow \infty$, since the inner sum consists of at most $N$ summands (see Assumption \ref{spec_1}). In conclusion, we have
 \begin{displaymath}
  D_p(x_{n_k},x_{n_l}) \rightarrow 0 \quad \text{for} \ l\rightarrow\infty,
 \end{displaymath}
 i.e., $\lbrace x_{n_k} \rbrace_{k \in \NN}$ is a Cauchy sequence and converges strongly. Its limit $\hat{x}:=\lim_{l\rightarrow \infty} x_{n_k}$ fulfills $\lVert F(\hat{x}) - y \rVert=0$ and thus $F(\hat{x})=y$, which implies $\hat{x} \in M_{F(x)=y}$. Since $x_{n_k} \rightharpoonup x^+$, we have $\hat{x}=x^+$. \\
 Now consider the sequence $\left\lbrace D_p(x^+,x_{n}) \right\rbrace_{n\in\NN}$. The continuity of $D_p(x^+,\cdot)$ yields
 \begin{displaymath}
  D_p(x^+,x_{n_k}) \rightarrow D_p(x^+,x^+) = 0 \quad \text{for} \ k \rightarrow \infty. 
 \end{displaymath}
 Consequently, the sequence $\left\lbrace D_p(x^+,x_{n}) \right\rbrace_{n\in\NN}$ which is convergent according to Proposition \ref{prop_descent_prop_algo1} has a subsequence converging to zero. Hence $\left\lbrace D_p(x^+,x_{n}) \right\rbrace_{n\in\NN}$ is also a null sequence, yielding the strong convergence of $\left\lbrace x_n \right\rbrace_{n\in\NN}$ by Theorem \ref{theorem_properties_bd} (v).
\end{proof}

\subsection{The RESESOP algorithm as a regularization method} \label{subsec_resesop_reg}

Let us now consider noisy data $y^{\delta}$ with noise level $\delta>0$. Again, we specify the search directions that are used in the respective algorithm. We now additionally require the Banach space $Y$ to be uniformly smooth, such that the duality mapping $J_2^Y$ is single-valued and continuous on bounded sets.

\vspace*{2ex}

\begin{assumption} \label{spec_2}
 \begin{itemize}
  \item[(i)] Let $n \in I_n^{\delta}$ for all iterations $n \in \NN$.
  \item[(ii)] Choose $I_n \subseteq \lbrace n-N+1,...,n \rbrace$ for some fixed $N>1$ such that the set of search directions
   \begin{displaymath}
    U^*_n := \left\lbrace u^*_{n,i} := F'(x_i)^*w_{n,i} \, : \, i \in I_n \right\rbrace
   \end{displaymath}
   is linearly independent.
   \item[(iii)] Let $I_n^{\delta} \subseteq \lbrace n-N+1,...,n \rbrace$ for some fixed $N>1$.
  \item[(iv)] Set $w_{n,n}^{\delta} := J^Y_2 \big( F\big(x_n^{\delta}\big) - y^{\delta} \big)$ for all $n\in \NN$.
 \end{itemize}
\end{assumption}

\vspace*{2ex}

\begin{lemma} \label{lemma_fsi}
 Algorithm \ref{algo2} stops after a finite number $n_* = n_*(\delta)$ of iterations, if the stopping criterion is the discrepancy principle with tolerance parameter
 \begin{equation} \label{cond_tau_ctc}
  \tau > \frac{1+\ctc}{1-\ctc},
 \end{equation}
 i.e., if
 \begin{displaymath}
  \lVert R_{n_*}^{\delta} \rVert \leq \tau \delta < \lVert R_{n}^{\delta} \rVert
 \end{displaymath}
 is fulfilled for all $n < n_*$.
\end{lemma}

\vspace*{2ex}

\begin{proof}
 As long as the discrepancy principle is not yet fulfilled at iteration $n\in\NN$, condition \eqref{cond_tau_ctc} implies
 \begin{displaymath}
  \lVert R_n^{\delta} \rVert > \tau\delta > \delta \cdot \frac{1+\ctc}{1-\ctc}
 \end{displaymath}
 and therefore
 \begin{displaymath}
 \begin{split}
  \alpha_{n,n}^{\delta} + \xi_{n,n}^{\delta} &= \left\langle u^*_{n,n}, x_n^{\delta} \right\rangle - \left\langle w_{n,n}^{\delta}, R_n^{\delta} \right\rangle + \big( \delta + \ctc (\lVert R_n^{\delta} \rVert + \delta) \big) \lVert w_{n,n}^{\delta} \rVert \\
  &= \left\langle u^*_{n,n}, x_n^{\delta} \right\rangle - \left\langle J_2^Y(R_{n}^{\delta}), R_n^{\delta} \right\rangle + \big( \delta + \ctc (\lVert R_n^{\delta} \rVert + \delta) \big) \lVert R_n^{\delta} \rVert \\
  &= \left\langle u^*_{n,n}, x_n^{\delta} \right\rangle - \left(\lVert R_n^{\delta}\rVert - \delta - \ctc (\lVert R_n^{\delta} \rVert + \delta) \right) \lVert R_n^{\delta} \rVert \\
  &< \left\langle u^*_{n,n}, x_n^{\delta} \right\rangle.
 \end{split}
 \end{displaymath}
 By projecting first onto the stripe $H\big(u^*_{n,n}, \alpha_{n,n}^{\delta}, \xi_{n,n}^{\delta}\big)$, i.e., onto the upper bounding hyperplane $H\big(u^*_{n,n}, \alpha_{n,n}^{\delta}+\xi_{n,n}^{\delta}\big)$, we obtain, analogously to the proof of Proposition \ref{prop_descent_prop_algo1}, the descent property
 \begin{equation} \label{desc_prop_reg}
  D_p(x_{n+1}^{\delta},z) \leq D_p(x_{n}^{\delta},z) - \frac{(1-\ctc - \tau^{-1}\cdot(1+\ctc))^p}{p G_{p^*}^{p-1}\cdot c_F^p} \cdot \lVert R_n^{\delta} \rVert^p
 \end{equation}
 for $z\in M_{F(x)=y}$. Let us now assume that there is no finite stopping index $n_* \in \NN$, such that $\lVert R_n^{\delta} \rVert > \tau\delta$ for all $n \in \NN$. From a calculation as in the proof of Proposition \ref{prop_weak_conv_algo1} we derive, however, that the sequence $\lbrace\lVert R_n^{\delta} \rVert \rbrace_{n\in\NN}$ is a null sequence, which is a contradiction to our assumption. As a consequence, there exists an $n_* \in \NN$ such that
 \begin{displaymath}
  \lVert R_{n_*}^{\delta} \rVert \leq \tau\delta
 \end{displaymath}
 and the discrepancy principle is fulfilled.
\end{proof}

\vspace*{2ex}

\begin{proposition}\label{prop_stab}
 The iterates $x_n^{\delta}$, $n\in\NN$, depend continuously on the data $y^{\delta}$, i.e.,
 \begin{equation}
  x_n^{\delta} \rightarrow x_n \quad \text{for} \ \delta \rightarrow 0,
 \end{equation}
 where $\lbrace x_n \rbrace_{n\in\NN}$ is the sequence of iterates generated by Algorithm \ref{algo1} for the respective exact data $y$.
\end{proposition}

\vspace*{2ex}

\begin{proof}
 Let $n \in \NN$ be a fixed iteration index. We prove by induction that $x_n^{\delta}$ depends continuously on $y^{\delta}$. In the first step for $k=0$ we have $I_0=I_0^{\delta} = \lbrace 0 \rbrace$ and 
 \begin{displaymath}
  x_1^{\delta} = J_{p^*}^{X^*}\left( J_p^X(x_0) - t_{0,0}^{\delta} F'(x_0)^* J_2^Y(F(x_0) - y^{\delta}) \right).
 \end{displaymath}
 According to Theorem \ref{theorem_properties_dm}, the occurring duality mappings are continuous since $X$ and $Y$ are uniformly smooth and $X^*$ is $p^*$-smooth. 
 Consequently, due to the continuity of the operator $F'(x_0)^*$, the search direction $u^*_{0,0}$ depends continuously on $y^{\delta}$. The optimization parameter $t_{0,0}^{\delta}$ is now calculated by minimizing the strictly convex smooth functional $h$ (see also Proposition \ref{prop_proj}), yielding its continuous dependence on $y^{\delta}$. Thus, the iterate $x_1^{\delta}$ fulfills $x_1^{\delta} \rightarrow x_1$ for $\delta \rightarrow 0$. \\
 We now assume that the iterates $x_i^{\delta}$, $0 \leq i \leq k$ depend continuously on $y^{\delta}$. By the same arguments as before, in addition to the continuity of the mapping $x \mapsto F'(x)$ and the linear independence of the search directions, we obtain the continuous dependence of
 \begin{displaymath}
  x^{\delta}_{n+1} = J_{p^*}^{X^*} \Big(J^X_p\big(x^{\delta}_n\big) - \sum_{i\in I^{\delta}_n} t^{\delta}_{n,i} F'\big(x^{\delta}_i\big)^*w^{\delta}_{n,i} \Big),
 \end{displaymath}
 on the data $y^{\delta}$ from the inductive hypothesis.
\end{proof}

\vspace*{2ex}

\begin{theorem}
 For given noisy data $y^{\delta}\in Y$ Algorithm \ref{algo2} with the specifications from Assumption \ref{spec_2} yields a regularized solution $x_{n_*(\delta)}^{\delta}$ of the nonlinear inverse problem \eqref{op_eq}, i.e., we have
 \begin{equation} \label{reg}
  x_{n_*(\delta)}^{\delta} \rightarrow x^+ \quad \text{for} \ \delta \rightarrow 0,
 \end{equation}
 if $x_0^{\delta} \in B^p_{\varrho}(x^+)$, where $x^+ \in M_{F(x)=y}$.
\end{theorem}

\vspace*{2ex}

\begin{proof}
 Let $y^{\delta}\rightarrow y$ for $\delta \rightarrow 0$, where $y$ are the exact data of the operator equation $F(x)=y$. We choose a null sequence $\lbrace \delta_j \rbrace_{j\in\NN} \subseteq \RR_{>0}$, fulfilling
 \begin{displaymath}
  \lVert y^{\delta_j} - y \rVert \leq \delta_j
 \end{displaymath}
 for all $j \in \NN$. For each noise level $\delta_j$, we define the respective finite stopping index
 \begin{displaymath}
  n_j := n_*(\delta_j),
 \end{displaymath}
 which exists according to Lemma \ref{lemma_fsi} and fulfills
 \begin{displaymath}
  \lVert R_{n_j}^{\delta_j} \rVert \leq \tau \delta_j.
 \end{displaymath}
 The validity of the descent property yields $\lbrace x_{n_j}^{\delta_j} \rbrace_{j\in\NN} \subseteq B^p_{\varrho}(x^+)$, which shows that the sequence $\lbrace x_{n_j}^{\delta_j} \rbrace_{j\in\NN}$ is bounded and has weak cluster points. Since we have 
 \begin{displaymath}
  \left\lVert F\big(x_{n_j}^{\delta_j}\big) - y^{\delta_j} \right\rVert = \left\lVert R_{n_j}^{\delta_j} \right\rVert \leq \tau\delta_j \rightarrow 0
 \end{displaymath}
 for $j \rightarrow \infty$, the weak sequential closedness of $F$ yields that the weak cluster points are contained in $M_{F(x)=y}$. \\
 We now assume, without loss of generality, that the sequence $\lbrace n_j \rbrace_{j\in\NN}$ is increasing. For every $k\in\NN$ we find $n_k, j_k \in \NN$ such that
 \begin{equation}
  D_p(x_{n_k}, x_{n_k}^{\delta_j}) \leq \frac{1}{k}, \quad D_p\big( x_{n_k}, x^+ \big) \leq \frac{1}{k}
 \end{equation}
 for all $j \geq \max\lbrace k, n_k \rbrace$, due to our stability result Proposition \ref{prop_stab}. Finally, we estimate
 \begin{displaymath}
  \begin{split}
   D_p\big( x_{n_j}^{\delta_j}, x^+ \big) &\leq D_p\big( x_{n_k}^{\delta_j}, x^+ \big) \\
                                        &\leq D_p\big( x_{n_k}, x^+ \big) - D_p\big( x_{n_k}, x_{n_k}^{\delta_j} \big) \\
                                        &\leq D_p\big( x_{n_k}, x^+ \big) + \frac{1}{k} \\
                                        &\leq \frac{2}{k}  \quad \rightarrow 0 \quad \text{for} \ k \rightarrow \infty,
  \end{split}
 \end{displaymath}
 using the convergence of the unperturbed sequence of iterates, see Theorem \ref{theorem_conv_exact_data}, and the descent property \eqref{desc_prop}. Theorem \ref{theorem_properties_bd} (v) finally yields the desired result \eqref{reg}.
\end{proof}

\vspace*{2ex}

\begin{remark}
 The Algorithms \ref{algo1} and \ref{algo2} yield a regulation of the step width for the Landweber method, if we choose $I_n = I_n^{\delta}=\lbrace n\rbrace$ for all $n\in\NN$. This method has been explicitly addressed by Maass and Strehlow \cite{pmrs16} for nonlinear inverse problems in Banach spaces in combination with sparsity constraints. 
\end{remark}

\section{A numerical example} \label{sec_numex}
We want to present a numerical evaluation of the proposed method. To this end we consider a fast algorithm with two search directions, which has been introduced before for linear operators in Banach spaces \cite{ss09} as well as for nonlinear problems in Hilbert spaces, see \cite{ws16}. Applications of this algorithm have shown a significant reduction of reconstruction time, in particular in terahertz tomography \cite{awts17}. We start by introducing this algorithm for nonlinear inverse problems in Banach spaces and apply it to a parameter identification problem, comparing the results to the Landweber method. 

\subsection{RESESOP with two search directions}
The following algorithm is a special case of Algorithm \ref{algo2} with the specifications
\begin{equation}\label{index_set_algo3}
 I_0^{\delta} = \lbrace 0 \rbrace, \quad I_n^{\delta} := \lbrace n-1,n \rbrace \ \text{for all} \ n\geq 1.
\end{equation}
The first iteration thus represents a Landweber step with a regulation of the step width, obtained by projecting the initial value onto its corresponding stripe. Note that Algorithm \ref{algo3} is presented as a regularization method for reconstructions from noisy data $y^{\delta}$. The corresponding method for exact data is obtained by setting $\delta=0$ and defining a suitable stopping rule. \\
According to Proposition 4 from \cite{ss09}, which is a general statement for Bregman projections in Banach spaces and is proved in a similar way as Theorem \ref{prop_descent_prop_algo1}, the Bregman projection of some $x\in X$ onto the intersection of two halfspaces (or stripes) can be computed by at most two Bregman projections onto (intersections of) the respective bounding hyperplanes, if $x$ is contained in one of the halfspaces (one of the stripes). This can be easily illustrated in a Hilbert space setting, see \cite{ss09}. 

\vspace*{2ex}

\begin{algorithm} \label{algo3}
 Use the definition \eqref{index_set_algo3} in Algorithm \ref{algo2} and set
 \begin{displaymath}
  w^{\delta}_{i} := w^{\delta}_{n,i} := J^Y_2\big(F(x_i^{\delta})-y^{\delta}\big) = J_2^{Y}\big(R_i^{\delta}\big)
 \end{displaymath}
 as well as $u^*_i := u^*_{n,i}$, $\alpha^{\delta}_i := \alpha^{\delta}_{n,i}$, $\xi^{\delta}_i := \xi^{\delta}_{n,i}$, and $H^{\delta}_{i}:=H^{\delta}_{n,i}$ for all $n\in\NN$ and $i\in I_n^{\delta}$. Choose an initial value $x_0^{\delta}:=x_0\in X$. As long as $\lVert R_n^{\delta} \rVert > \tau\delta$, we have
 \begin{displaymath}
  x_n^{\delta} \in H_>\big(u_n^*,\alpha^{\delta}_n + \xi^{\delta}_n\big) \cap H_{n-1}^{\delta}. 
 \end{displaymath}
 Now compute the new iterate $x_{n+1}^{\delta}$ according to the following two steps.
 \begin{itemize}
  \item[(i)] Compute
  \begin{displaymath}
   \tilde{x}_{n+1}^{\delta} := \Pi^p_{H(u_n^*,\alpha^{\delta}_n + \xi^{\delta}_n)}\big(x_n^{\delta}\big),
  \end{displaymath}
  i.e., compute $J_p(\tilde{x}_{n+1}^{\delta}) = J_p(x_n^{\delta}) - t_n^{\delta} u_n^*$, where $t_n^{\delta}$ minimizes 
  \begin{displaymath}
   h_{n,1}(t) := \frac{1}{p^*} \left\lVert J_p(x_n^{\delta}) - t_n^{\delta} u_n^* \right\rVert^{p^*} + t \cdot \left( \alpha_n^{\delta} + \xi_n^{\delta} \right).
  \end{displaymath}
  The intermediate iterate $\tilde{x}_{n+1}^{\delta}$ then fulfills the descent property
  \begin{equation} \label{desc_prop_algo3_i}
   D_p\big(\tilde{x}^{\delta}_{n+1},z\big) \leq D_p\big(x^{\delta}_n,z\big) - \frac{\big(\lVert R_n^{\delta} \rVert - \delta - \ctc (\lVert R_n^{\delta} \rVert + \delta)\big)^p}{p G_{p^*}^{p-1}\cdot c_F^p} 
  \end{equation}

  If $\tilde{x}_{n+1}^{\delta} \in H_{n-1}^{\delta}$, we are done and set
  $x_{n+1}^{\delta} := \tilde{x}_{n+1}^{\delta}$. Otherwise, we proceed with step (ii).
  \item[(ii)] First we decide whether $\tilde{x}_{n+1}^{\delta} \in H_>\big(u^*_{n-1}, \alpha_{n-1}^{\delta}+\xi_{n-1}^{\delta}\big)$ \\ or $\tilde{x}_{n+1}^{\delta} \in H_<\big(u^*_{n-1}, \alpha_{n-1}^{\delta}-\xi_{n-1}^{\delta}\big)$. Calculate accordingly
  \begin{displaymath}
   x_{n+1}^{\delta} := \Pi^p_{H(u_n^*,\alpha^{\delta}_n + \xi^{\delta}_n) \cap H(u^*_{n-1}, \alpha_{n-1}^{\delta}\pm\xi_{n-1}^{\delta})}\big( \tilde{x}_{n+1}^{\delta} \big),
  \end{displaymath}
  i.e., determine
  \begin{displaymath}
   J_p(x_{n+1}^{\delta}) = J_p(\tilde{x}_{n+1}^{\delta}) - t_{n,n}^{\delta} u_n^* - t_{n,n-1}^{\delta} u_{n-1}^*,
  \end{displaymath}
  where $\big( t_{n,n}^{\delta}, t_{n,n-1}^{\delta} \big)$ minimizes
  \begin{displaymath}
   h_{n,2}(t) := \frac{1}{p^*} \left\lVert J_p(x_n^{\delta}) - t_1^{\delta} u_n^* - t_2^{\delta} u_{n-1}^* \right\rVert^{p^*} + t_1 \cdot \left( \alpha_n^{\delta} + \xi_n^{\delta} \right) + t_2 \cdot \left( \alpha_{n-1}^{\delta} \pm \xi_{n-1}^{\delta} \right).
  \end{displaymath}
  This yields
  \begin{displaymath}
   x_{n+1}^{\delta} = \Pi^p_{H_n^{\delta}\cap H_{n-1}^{\delta}}\big(x_n^{\delta}\big)
  \end{displaymath}
  and the descent property
  \begin{equation} \label{desc_prop_algo3_ii}
   D_p\big(x^{\delta}_{n+1},z\big) \leq D_p\big(x^{\delta}_n,z\big) - \frac{1}{pG_{p^*}^{p-1}} S_n^{\delta},
  \end{equation}
  where
  \begin{equation} \label{sn_algo3}
  \begin{split}
   S_n^{\delta} &:= \left(\frac{\big(\lVert R_n^{\delta} \rVert - \delta - \ctc (\lVert R_n^{\delta} \rVert + \delta)\big)}{c_F}\right)^p \\
                &\qquad + \left( \frac{\left\lvert \left\langle u_{n-1}^*, \tilde{x}_{n+1}^{\delta} \right\rangle - \big( \alpha_{n-1}^{\delta} \pm \xi_{n-1}^{\delta} \big)\right\rvert}{\gamma^{\delta}_n \lVert u_{n-1}^* \rVert} \right)^p
  \end{split}
  \end{equation}
  and
  \begin{displaymath}
   \gamma_n^{\delta} := \left( 1 - \frac{1}{(p-1) G_{p^*}^{p-1}} \left( \frac{\left\lvert \left\langle u_n^*, J_{p^*}^{X^*} (u_{n-1}^*) \right\rangle \right\rvert}{\lVert u_n^* \rVert \cdot \lVert J_{p^*}^{X^*} (u_{n-1}^*) \rVert } \right)^p \right)^{\frac{1}{p^*}}.
  \end{displaymath}
 \end{itemize}
\end{algorithm}

 \vspace*{2ex}

\begin{remark}
 \begin{itemize}
  \item[(i)] The statements in Algorithm \ref{algo3}, in particular the descent properties \eqref{desc_prop_algo3_i} and \eqref{desc_prop_algo3_ii}, are a direct consequence of Proposition 4 from \cite{ss09}, which is a general statement for the Bregman projection of $x \in X$ onto the intersection of two halfspaces. In the case of nonlinear forward operators, one proceeds as in the proof of Proposition \ref{prop_descent_prop_algo1}.
  \item[(ii)] The estimate for the descent property \eqref{desc_prop_algo3_ii} indicates the advantage of using multiple search directions: the second term in \eqref{sn_algo3} quantifies the acceleration due to the second search direction.
   \item[(iii)] Algorithm \ref{algo3} has previously been formulated for linear operators in Banach spaces \cite{ss09} and for nonlinear operators in Hilbert spaces \cite{ws16}. These algorithms are based on the observation that the Bregman projection of $x\in X$ onto the intersection of two halfspaces can be uniquely determined by at most two projections onto (intersections of) the bounding hyperplanes if $x$ is already contained in one of the halfspaces.
 \end{itemize}
\end{remark}

\subsection{Parameter identification with SESOP/RESESOP}
We evaluate the proposed method by solving a well-known nonlinear parameter identification problem, using Algorithm \ref{algo3}. The performance is compared to the Landweber-type method that is obtained from Algorithms \ref{algo1} and \ref{algo2} with a single search direction (the Landweber direction) per iteration. \\
Consider the boundary value problem
\begin{equation} \label{bvp}
\begin{aligned}
 - \laplace u + cu &= f &&\quad \text{in} \ \Omega, \\
                 u &= g &&\quad \text{on} \ \partial\Omega
\end{aligned}
\end{equation}
for known functions $f \in L^{2}(\Omega)$ and $g \in H^{\frac{3}{2}}(\partial\Omega)$. We want to reconstruct the parameter $c \in L^{2}(\Omega)$ from the knowledge of the (possibly perturbed) solution $u \in H^2(\Omega)$. The nonlinear forward problem is formulated as
\begin{displaymath}
 F: \mathcal{D}(F) \subseteq X:=L^2(\Omega) \rightarrow Y:=H^2(\Omega), \ F(c) = u,
\end{displaymath}
where $u$ solves \eqref{bvp} and
\begin{displaymath}
 \mathcal{D}(F) = \big\lbrace c \in L^2(\Omega) \, : \, \lVert c-\hat{c} \rVert_{L^2(\Omega)} \leq \gamma \ \text{for a} \ \hat{c} \in L^{\infty}(\Omega) \ \text{with} \ \hat{c} \geq 0 \ \text{a.e.}  \big\rbrace
\end{displaymath}

This parameter identification problem is a standard problem commonly used to evaluate reconstruction techniques for nonlinear inverse problems. In particular, the Fr\'echet derivative and its adjoint are known, see, e.g., \cite{hns_cl, kss09, kaltenbacher97, skhk12} for an analysis of this boundary value problem, also in the context of Banach spaces.

The Fr\'echet derivate is given by
\begin{equation}
 F'(c): L^2(\Omega) \rightarrow H^2(\Omega), \ F'(c)h = - L(c)^{-1} (h\cdot F(c)),
\end{equation}
where
\begin{displaymath}
 L(c) : H^2(\Omega) \rightarrow L^2(\Omega), \ L(c)u = - \laplace u + cu.
\end{displaymath}
Its dual is expressed as
\begin{displaymath}
 F'(c)^*: L^{2}(\Omega) \subseteq (H^2(\Omega))^* \rightarrow L^{2}(\Omega), \ F'(c)^* w = - u \cdot L(c)^{-1} w
\end{displaymath}
for $w \in L^{2}(\Omega)$.

\vspace*{2ex}

Instead of the Hilbert spaces $L^2(\Omega)$ and $H^2(\Omega)$, we now want to consider Banach spaces as source and data space. We use $u \in H^2(\Omega) \subseteq W^{2,s}(\Omega)$, $1 < s \leq 2$, as well as $c \in L^r(\Omega) \subseteq L^2(\Omega)$, $2 \leq r < \infty$ and consider the operator
\begin{displaymath}
 F : \mathcal{D}(F) \subseteq X := L^r(\Omega) \rightarrow Y := L^s(\Omega)
\end{displaymath}
as a mapping between Banach spaces $X$ and $Y$, see also \cite{kss09}.

\vspace*{2ex}

Note that according to Example \ref{ex1}, the spaces $L^r(\Omega)$, $L^s(\Omega)$, $L^{r^*}(\Omega)$, and $L^{s^*}(\Omega)$ are uniformly smooth and uniformly convex, such that the duality mappings are single-valued.

\vspace*{2ex}

For the implementation, we choose $\Omega = (0,1)\times(0,1)$ and use finite differences for the discretization, see, e.g., \cite{wh17}. The resulting grid $\Omega_h$ has $N+1$ grid points in $x$- as well as in $y$-direction with step size $h = (N+1)^{-1}$. The discretization of a function $f$ on $\Omega$ is denoted by $f^h$. In order to approximate the $L^p$-norm on $\Omega_h$, we use the \emph{weighted $p$-norm}
\begin{displaymath}
 \left\lVert f^h \right\rVert_{p,h} = h^{\frac{2}{p}} \cdot \left( \sum_{i,j=0}^{N+1} \lvert f(ih,jh) \rvert^p \right)^{\frac{1}{p}}.
\end{displaymath}

In order to have a precise representation of the exact solution $u$ and the parameter $c$ of \eqref{bvp}, we provide $u$ and $c$ analytically, such that $f$ and the boundary values $g$ are uniquely determined by
\begin{equation}
 f = - \laplace u + cu, \quad g = u |_{\partial\Omega}.
\end{equation}
Let 
\begin{displaymath}
 u \, : \, \overline{\Omega} \rightarrow \RR, \ u(x,y) = 16 x(x-1)y(1-y) + 1
\end{displaymath}
and
\begin{displaymath}
 c \, : \, \overline{\Omega} \rightarrow \RR, \  c(x,y) = \frac{3}{2}\sin(2\pi x)\sin(3\pi y) + 3\left( \left( x - \frac{1}{2} \right)^2 + \left( y - \frac{1}{2} \right)^2 \right) + 2.
\end{displaymath}
In our iterations, we will use the function
\begin{displaymath}
 c_0 \, : \, \overline{\Omega} \rightarrow \RR, \  c_0(x,y) = 3\big((x-0.5)^2 + (y-0.5)^2\big) + 2 + 8x(x-1)y(1-y)
\end{displaymath}
as an initial value.

\begin{figure}[ht]
 \centering
 \subfloat[]{\includegraphics[width=0.45\textwidth]{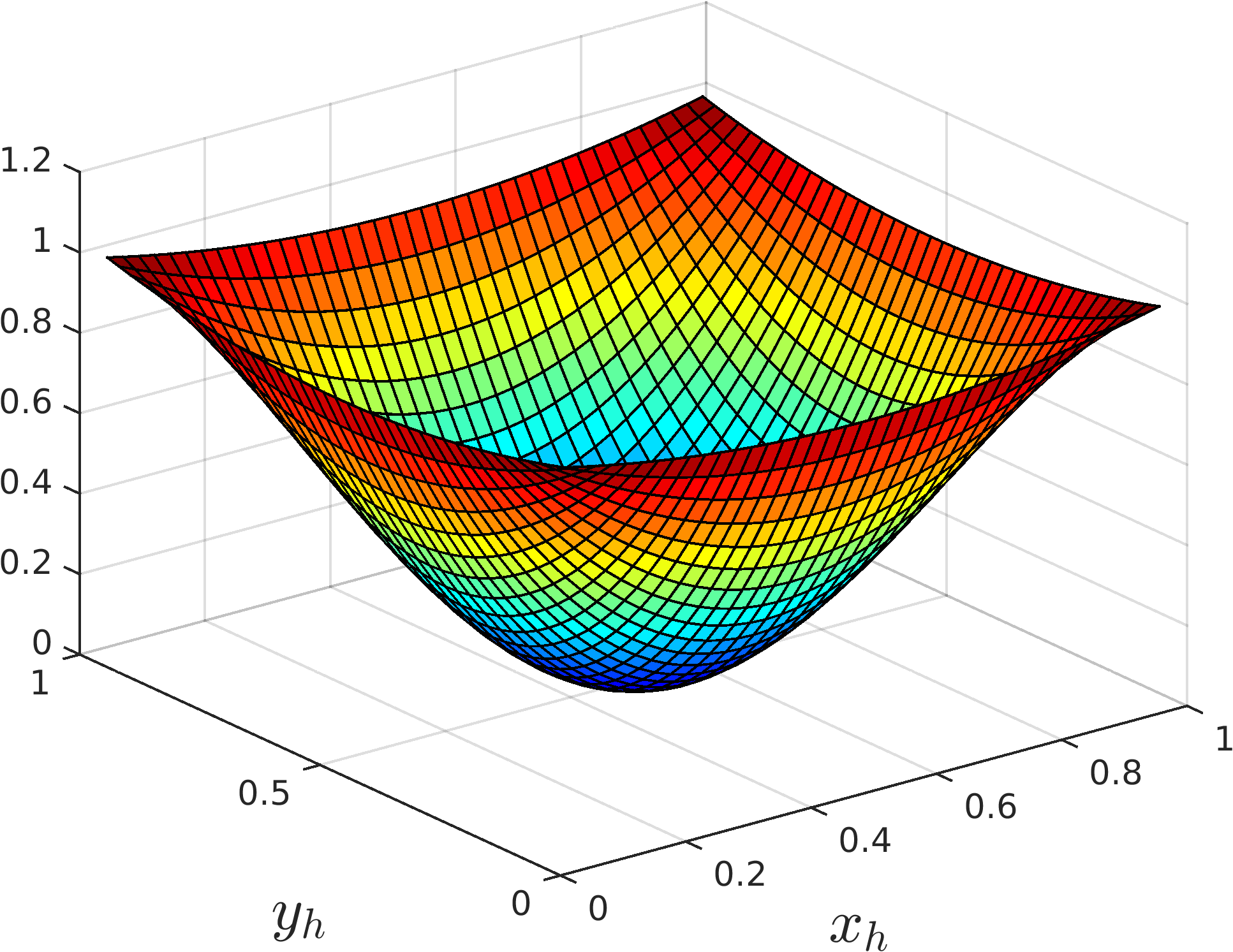}} \hspace*{2ex}
 \subfloat[]{\includegraphics[width=0.45\textwidth]{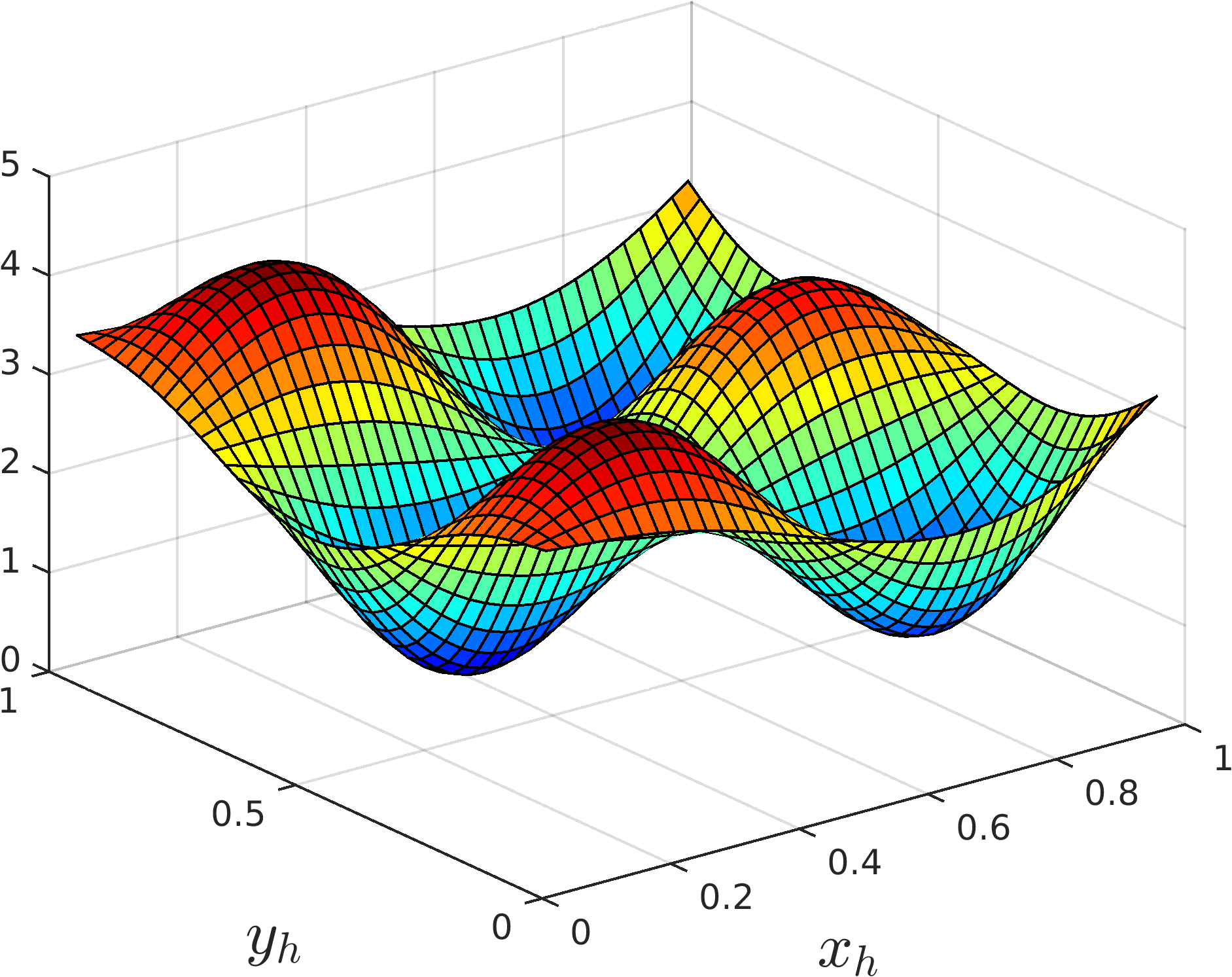}}
 \caption{\label{figure_u_c}Surface plot of (a) exact function $u_h$ and (b) exact parameter $c_h$.}
\end{figure}

For our reconstructions, we use synthetic exact data $u^h$ and synthetic noisy data $u^{h,\delta}$, which is obtained by adding a noise matrix $v$ to the exact data matrix $u_h$. To this end, we generate a matrix $v$ with random, equally distributed entries in $[-1,1]$ and set
\begin{equation}
 u_h^{\delta} = u_h + \delta \cdot \frac{v}{\lVert v \rVert_{p,h}}, 
\end{equation}
which ensures that $\left\lVert u_h^{\delta} - u_h \right\rVert_{p,h} = \delta$.

\vspace*{1ex}

The reconstructions are performed iteratively by
\begin{displaymath}
 c_{n+1}^{h,\delta} = c_n^{h,\delta} - d_n^{h,\delta},
\end{displaymath}
where $n$ is the iteration index, $\left\lbrace c_n^{h,\delta} \right\rbrace_{n\in\NN}$ is the sequence of iterates and $d_n^{h,\delta}$ is calculated according to
\begin{itemize}
 \item[(A)] the Landweber-type iteration that is based on the successive projection onto stripes, i.e., Algorithm \ref{algo1} or \ref{algo2} with $I_n=\lbrace n \rbrace$ and a single search direction (the current gradient $g_n^{\delta}$, see \eqref{current_gradient}) per iteration $n$, and
 \item[(B)] the SESOP method \ref{algo3} with two search directions (the current gradient and the gradient from the previous step).
\end{itemize}
In the case of exact data, we omit the index $\delta$ in the above notation and stop the iteration as soon as the norm of the residual falls below a tolerance $T_Y>0$. For the generation of synthetic data, we use $N=50$, whereas for the reconstruction we choose $N=40$. In each iteration $n$, we calculate the \emph{relative error}
\begin{displaymath}
 e_{\mathrm{rel},n} := \frac{\lVert c^h_n - c^h \rVert_{p,h}}{\lVert c^h \rVert_{p,h}}
\end{displaymath}
in the reconstruction $c^h_n$.

\vspace*{2ex}

For our numerical experiment, we choose
\begin{equation}
 r=1.5, \ s=5, \quad \ctc = 0.01, \quad \tau = 1.1 \cdot \frac{1+\ctc}{1-\ctc}.
\end{equation}
The duality mappings $J_p(f)$ and $J_q(f)$, $1 < p,q < \infty$, in $f \in \RR^N$ are evaluated according to
\begin{displaymath}
 J_p(f)(x) = \lvert f(x) \rvert^{p-1} \mathrm{sign}(f(x)), \quad J_q(f)(x) = \lVert f(x) \rVert^{r-p} J_p(f)(x),
\end{displaymath}
see also \cite{ssl08}, Example 2.2.

\paragraph{Reconstructions from exact data}

We compare the performance of the methods (A) and (B) for exact data ($\delta=0$). The iteration is stopped at iteration $n_*$ if 
\begin{equation}
 \lVert R_{n_*} \rVert \leq T_Y := 5 \cdot 10^{-4} < \lVert R_{n} \rVert
\end{equation}
for all $n < n_*$, i.e., as soon as the norm of the residual falls below the bound $T_Y$ for the first time.

\begin{figure}[!ht]
 \centering
\subfloat[Reconstruction of $c_h$, method (A)]{
  \includegraphics[width=0.4\textwidth]{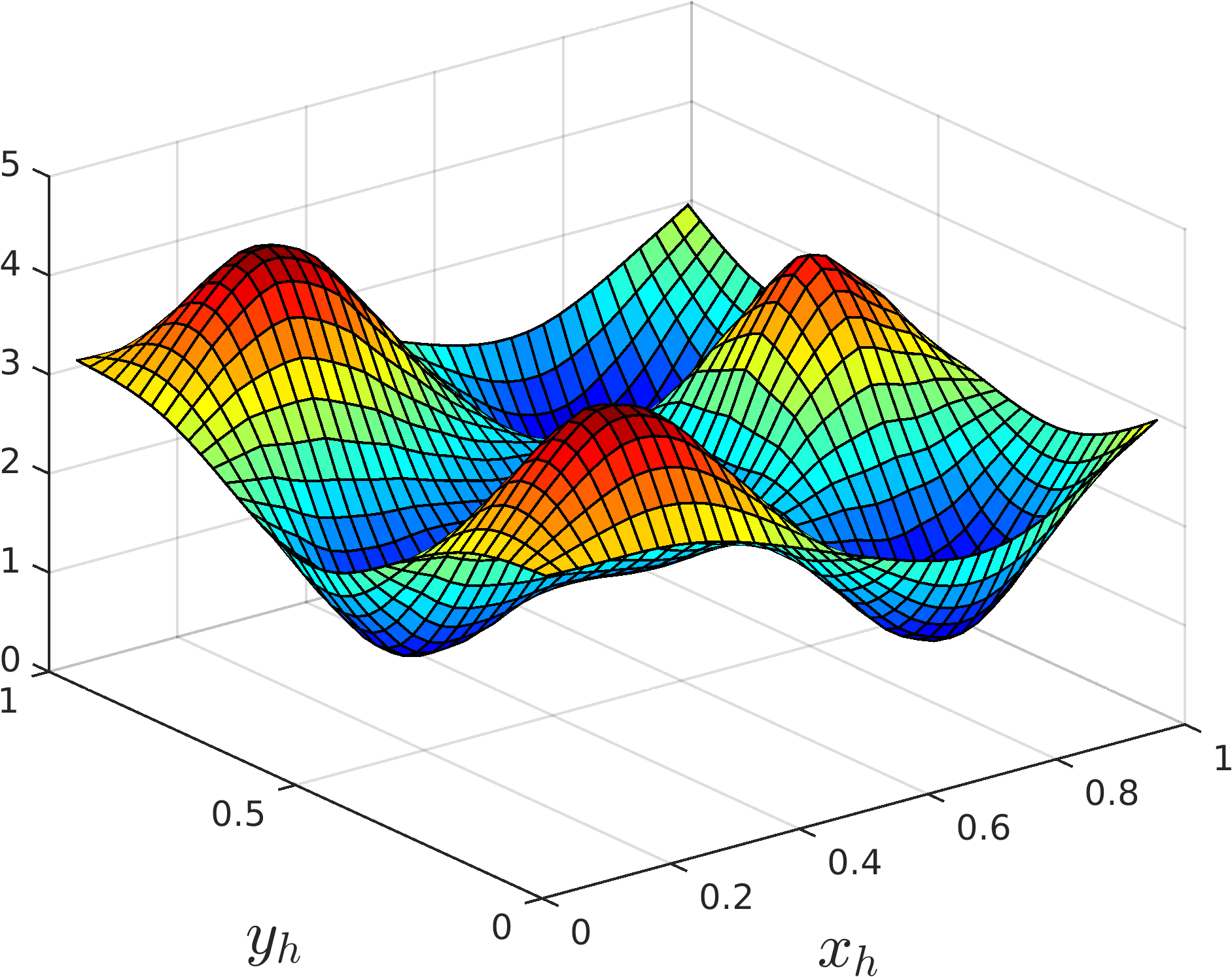}
}
\subfloat[Reconstruction of $c_h$, method (B)]{
  \includegraphics[width=0.4\textwidth]{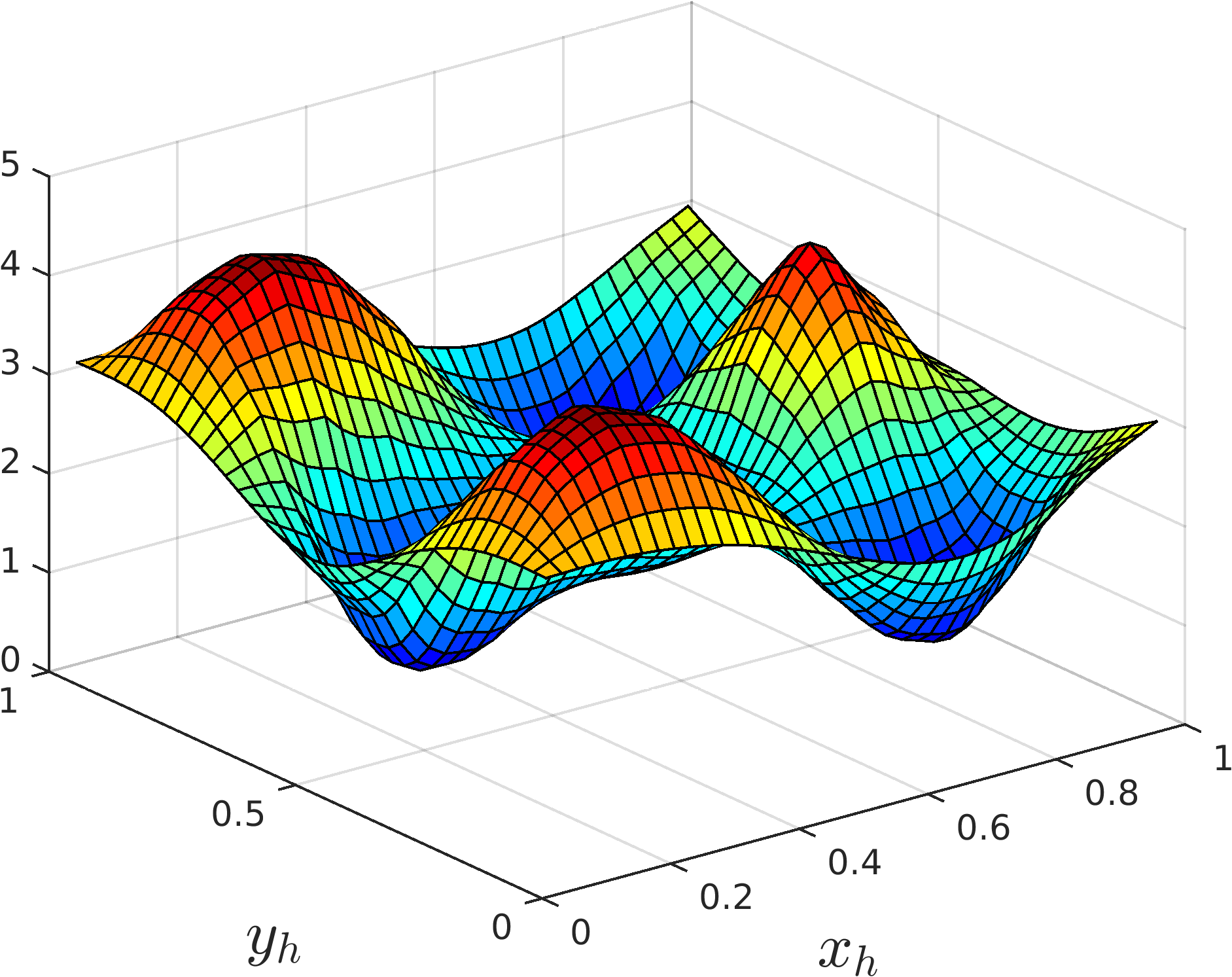}
}
\caption{\label{figure_reconstructions_toy_problem_exact} Reconstruction with the Landweber type method (A) and the RESESOP method (B): surface plots of the respective reconstructed parameter $c^h_{n_*}$.}
\end{figure}

\begin{table}[ht]
 \centering
 \begin{tabular}{|c|c|c|c|} 
  \hline
  Method  &  Number of iterations $n_*$ & Execution time      & Relative error $e_{\mathrm{rel},n_*}$ \\
  \hline
  (A)     &  $27$                       & $5.64 \mathrm{s}$   & $8.43 \%$ \\
  \hline
  (B)     &  $15$                       & $6.3 \mathrm{s}$    & $7.90 \%$ \\
  \hline
 \end{tabular}
 \caption{\label{table_comp_exact_BS}Some key data to evaluate the performance of the methods (A), (B) in the case of exact data $u_h$.}
\end{table}

Table \ref{table_comp_exact_BS} illustrates the performance of the two methods in this Banach space scenario. The respective reconstructions are displayed in Figure \ref{figure_reconstructions_toy_problem_exact}.

\begin{figure}[!ht]
\centering
\subfloat[Norm of residual versus iteration index, method (A)]{
  \includegraphics[width=0.4\textwidth]{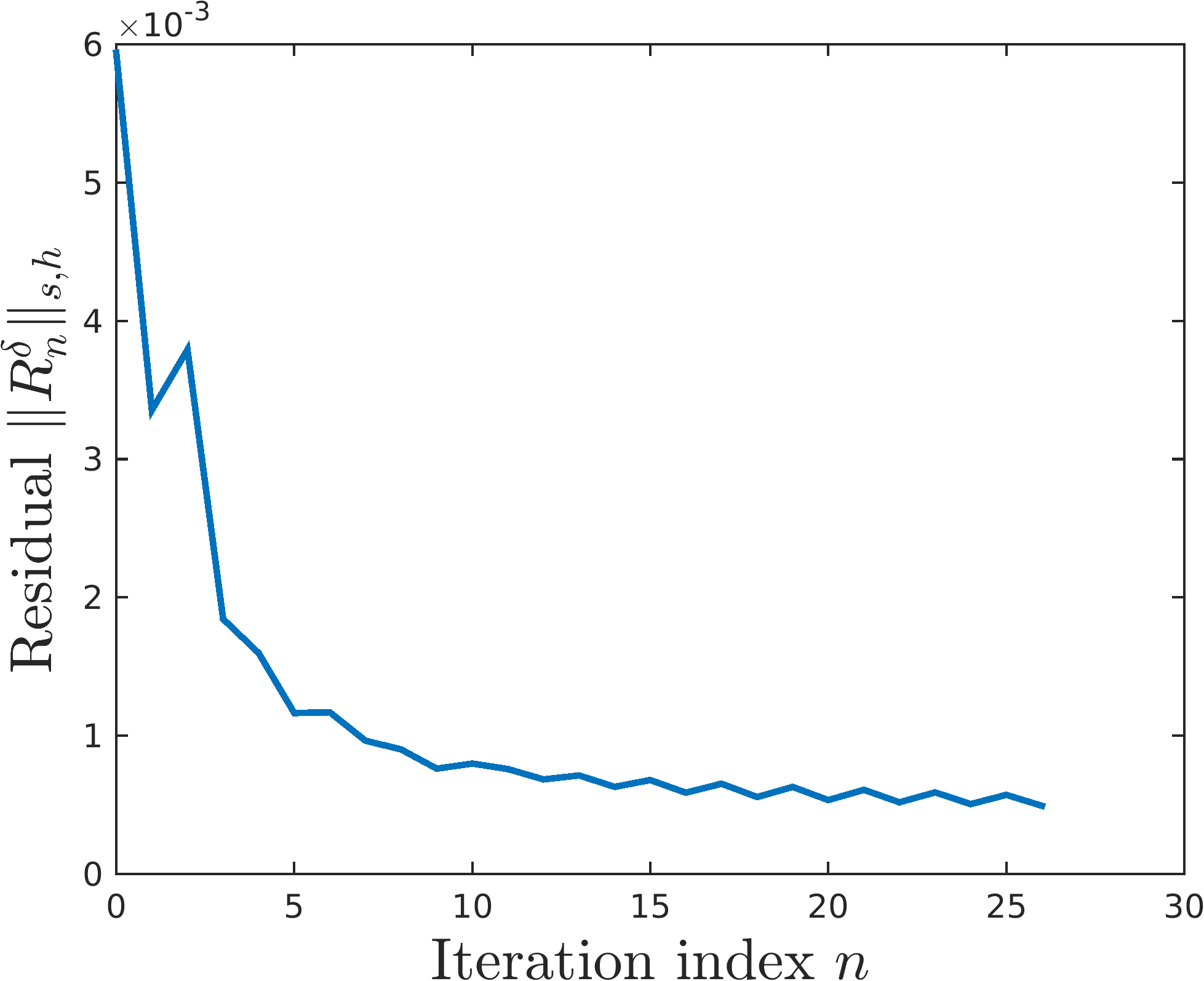}
} \hspace*{2ex}
\subfloat[Norm of residual versus iteration index, method (B)]{
  \includegraphics[width=0.4\textwidth]{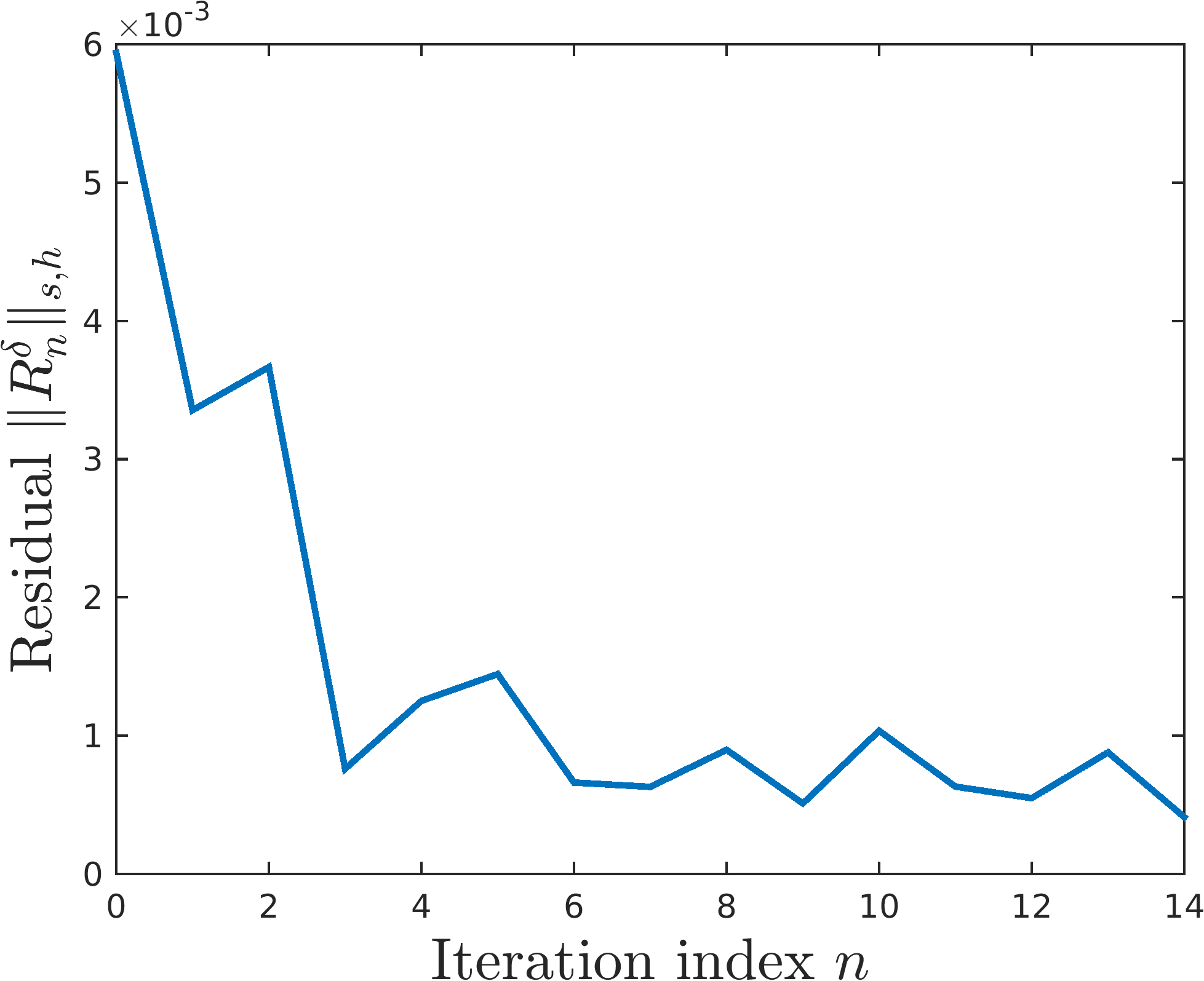}
}
\hspace{0mm}
\subfloat[Error versus iteration index, method (A)]{
  \includegraphics[width=0.4\textwidth]{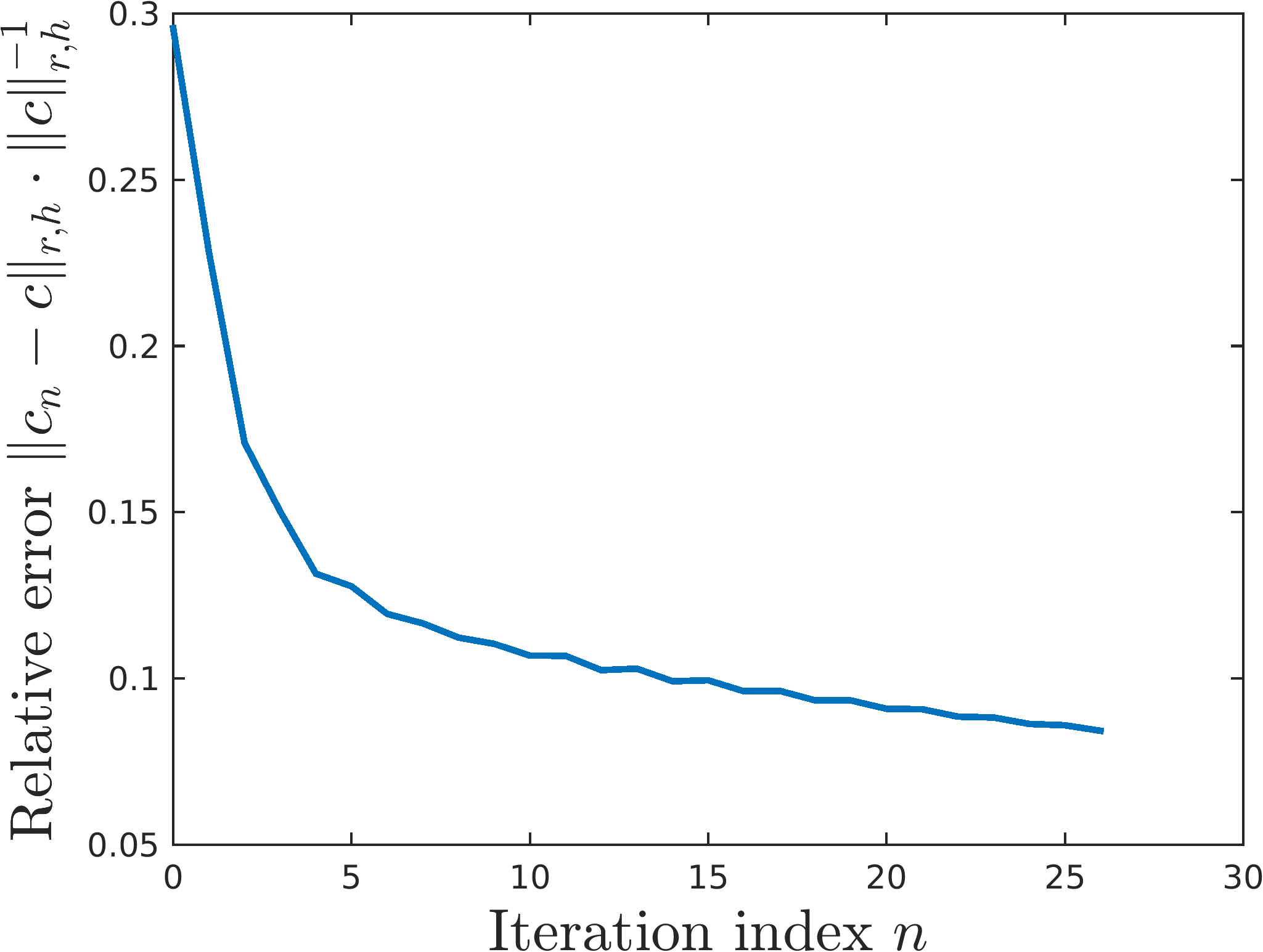}
}\hspace*{2ex}
\subfloat[Error versus iteration index, method (B)]{
  \includegraphics[width=0.4\textwidth]{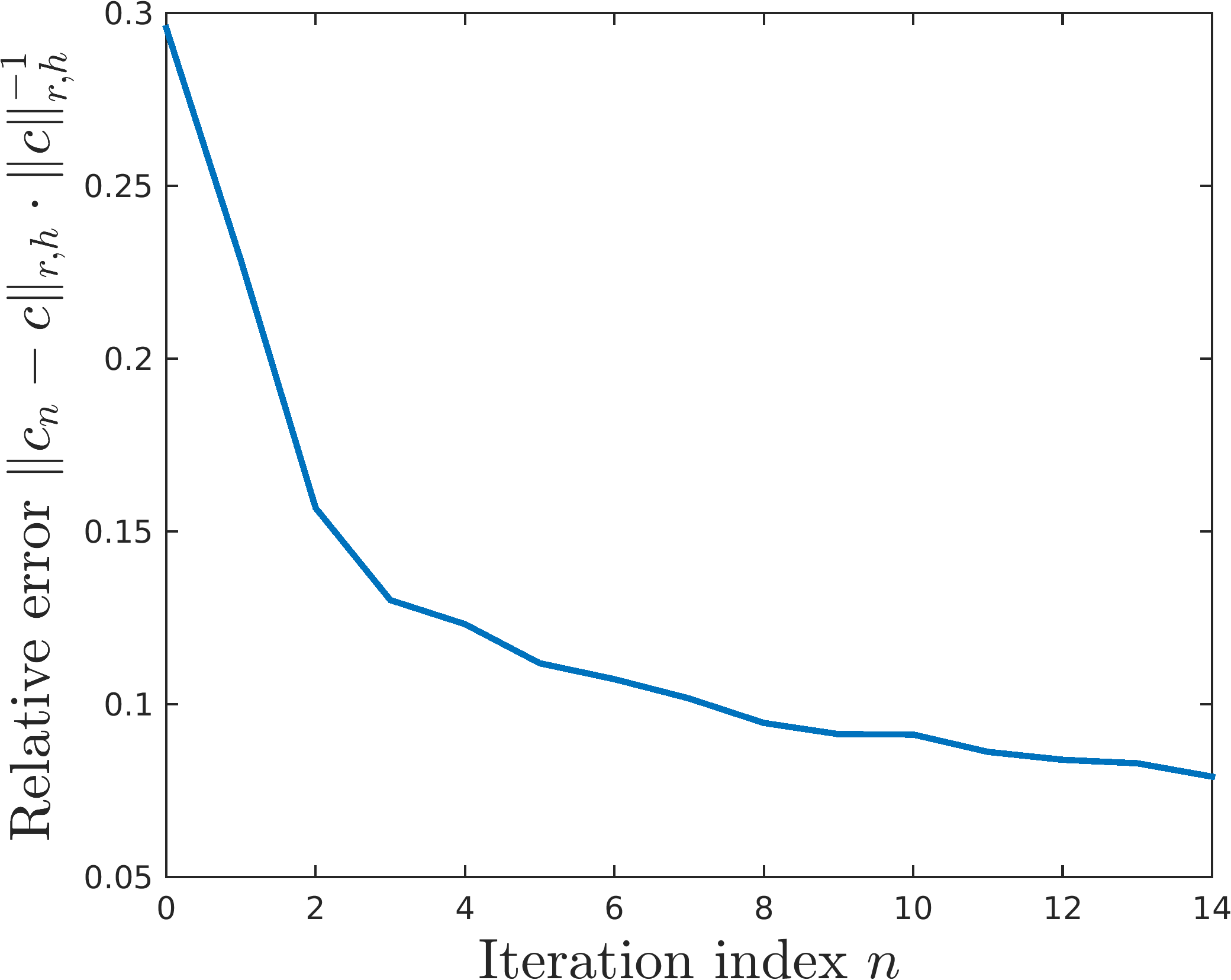}
}
\caption{\label{figure_residual_error}Comparison of residual and reconstruction error.}
\end{figure}

\vspace*{2ex}

\paragraph{Reconstructions from noisy data}

We now compare the methods (A) and (B) for noisy data with noise level $\delta = 5\cdot 10^{-4}$. The iteration is stopped by the discrepancy principle at iteration $n_*$. The performance of the methods is shown in Table \ref{table_comp_noisy_BS}.

\begin{table}[ht]
 \centering
 \begin{tabular}{|c|c|c|c|} 
  \hline
  Method  &  Number of iterations $n_*$ & Execution time      & Relative error $e_{\mathrm{rel},n_*}$ \\
  \hline
  (A)     &  $22$                       & $4.56 \mathrm{s}$   & $9.70 \%$ \\
  \hline
  (B)     &  $12$                       & $2.06 \mathrm{s}$    & $10.19 \%$ \\
  \hline
 \end{tabular}
 \caption{\label{table_comp_noisy_BS}Some key data to evaluate the performance of the methods (A), (B) in the case of exact data $u_h$.}
\end{table}

\begin{figure}[!ht]
 \centering
\subfloat[Reconstruction of $c_h$, method (A)]{
  \includegraphics[width=0.4\textwidth]{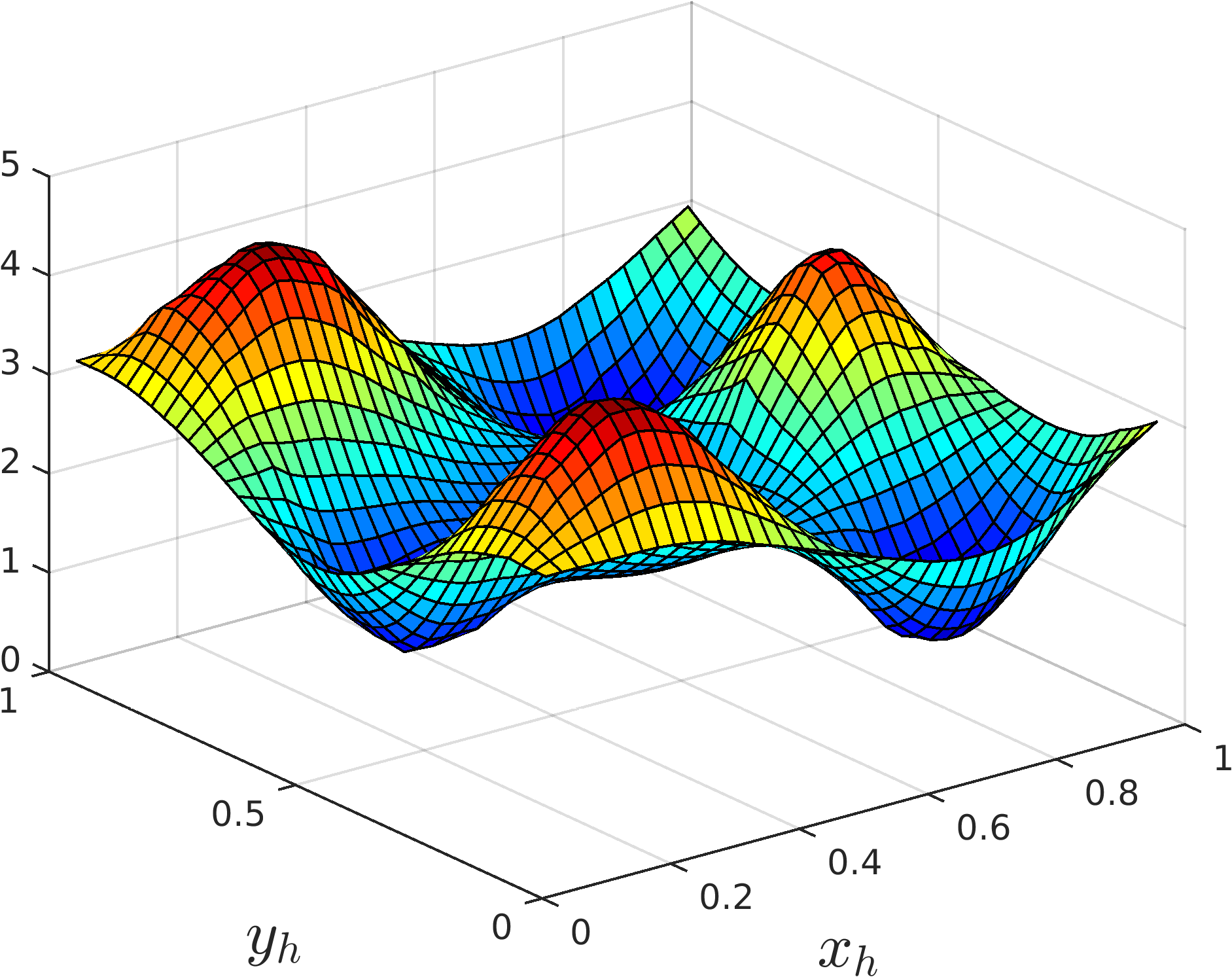}
}
\subfloat[Reconstruction of $c_h$, method (B)]{
  \includegraphics[width=0.4\textwidth]{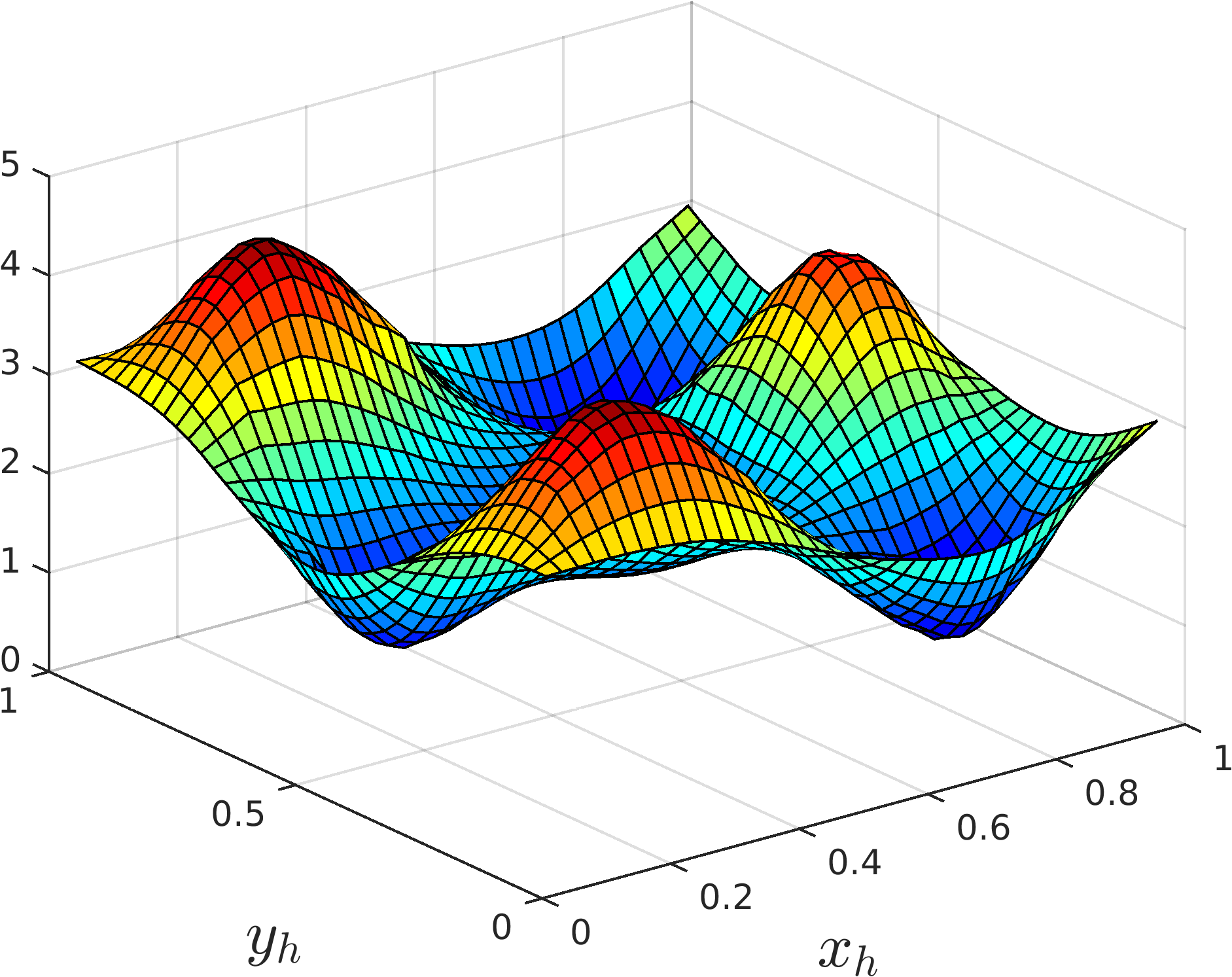}
}
\caption{\label{figure_reconstructions_toy_problem_noisy} Reconstruction with the Landweber type method (A) and the RESESOP method (B): surface plots of the respective reconstructed parameter $c^h_{n_*}$.}
\end{figure}

It is noteworthy that in step (ii) of Algorithm \ref{algo3}, the intermediate iterate $\tilde{x}_{n+1}^{\delta}$ is always contained in the halfspace above the stripe that belongs to the previous iterate. In Figure \ref{figure_residual_error} we see that the (relative) approximation error is monotonically decreasing during the iteration, indicating the validity of the descent property, whereas this is not the case for the norm of the residual. In addition, we observe the same effect of the choice of the norms as it has been described in \cite{ss09, skhk12}: If we set $r < 2$, $s > 2$, we reduce the total number of iterations in comparison to the case where $r>2$, $s<2$.

\section{Conclusion and outlook} \label{sec_co}

This work is dedicated to the adaption of the well-known sequential subspace optimization techniques to nonlinear inverse problems in Banach spaces. We have introduced a range of tools that are useful in the context of inverse problems in Banach spaces, including duality mappings, Bregman distances, and Bregman projections. The methods we discussed belong to the class of iterative techniques, and the new iterate is sought in the direction of a finite number of search directions. In each iteration, the set of search directions contains the current Landweber direction as well as Landweber directions from previous iterations. The iterates are chosen such that the current iterate is projected onto the (intersection of) locally defined stripes, whose width is determined by the noise level and the nonlinearity (i.e., the constant $\ctc$ from the tangential cone condition) of the forward operator. The proposed algorithm with two search directions has been successfully applied in parameter identification. \\
Possible extensions of this work include a more thorough analysis of search directions. For example, it might be fruitful to include directions $x_{m}^{\delta}-x_k^{\delta}$, $k < m \leq n$ in the search space as it has been done for linear inverse problems \cite{skhk12}. Generally, our method is designed for the numerical solution of nonlinear inverse problems, where the calculation of the Landweber direction is computationally expensive, such that a reduction of iterations yields a significant decrease in the time that is needed for the reconstruction. An example is the identification of the stored energy function from measurements of the displacement field in the context of structural health monitoring, which involves numerical evaluations of the elastic wave equations in each step of the Landweber iteration, see \cite{jsts16}.

\vspace*{2ex}

\textbf{Acknowledgement:} The research presented in this article was funded under the grant 05M16TSA by the German Federal Ministry of Education and Research (Bundesministerium f\"ur Bildung und Forschung, BMBF). 

\bibliographystyle{siam}


\pagestyle{myheadings}
\thispagestyle{plain}
\markboth{A. Wald}{Sequential subspace optimization for nonlinear inverse problems in Banach spaces}

\end{document}